\documentclass{amsart}
\usepackage{amssymb,amsmath,latexsym,epsfig}

\newtheorem{thm}{Theorem}[section]
\newtheorem{dfn}[thm]{Definition}
\newtheorem{cor}[thm]{Corollary}

\newtheorem{lemma}[thm]{Lemma}

\newcommand{\del}{\backslash}
\newcommand{\cl}{\hbox{\rm cl}}

\title[Lattice Path Matroids]{Lattice Path Matroids: Structural
Properties}
\date{\today}
\author[Joseph E.~Bonin]
       {Joseph E.~Bonin}
\address[Joseph E.~Bonin]
{Department of Mathematics\\ The George Washington University\\
Washington, D.C. 20052, USA} \email [Joseph E.~Bonin] {jbonin@gwu.edu}
\author[Anna de Mier]
{Anna de Mier}
\address[Anna de Mier]
{Mathematical Institute\\  24--29 St.~Giles \\ 
Oxford OX1 3LB, United Kingdom} \email [Anna de Mier] {demier@maths.ox.ac.uk}
                           
\begin{document}

\begin{abstract}
  This paper studies structural aspects of lattice path matroids, a
  class of transversal matroids that is closed under taking minors and
  duals.  Among the basic topics treated are direct sums, duals,
  minors, circuits, and connected flats.  One of the main results is a
  characterization of lattice path matroids in terms of fundamental
  flats, which are special connected flats from which one can recover
  the paths that define the matroid.  We examine some aspects related
  to key topics in the literature of transversal matroids and we
  determine the connectivity of lattice path matroids.  We also
  introduce notch matroids, a minor-closed, dual-closed subclass of
  lattice path matroids, and we find their excluded minors.
\end{abstract}

\maketitle

\section{Introduction}\label{sec:intro}

A lattice path matroid is a special type of transversal matroid whose
bases can be thought of as lattice paths in the region of the plane
delimited by two fixed bounding paths.  These matroids, which were
introduced and studied from an enumerative perspective in~\cite{lpm1},
have many attractive structural properties that are not shared by
arbitrary transversal matroids; this paper focuses on such properties.

The definition of lattice path matroids is reviewed in
Section~\ref{sec:back}, where we also give some elementary properties
of their bases and make some remarks on connectivity and
automorphisms. Section~\ref{sec:lpm} proves basic results that are
used throughout the paper; for example, we show that the class of
lattice path matroids is closed under minors, duals, and direct sums,
we determine which lattice path matroids are connected, and we
describe circuits and connected flats. The next section discusses
generalized Catalan matroids, a minor-closed, dual-closed subclass of
lattice path matroids that has particularly simple characterizations.
Section~\ref{sec:fundamental} introduces special connected flats
called fundamental flats that we use to characterize lattice path
matroids and to show that the bounding paths can be recovered from the
matroid.  In Section~\ref{sec:trans}, we describe the maximal
presentation of a lattice path matroid, and we use this result to give
a geometric description of these matroids as well as a polynomial-time
algorithm for recognizing lattice path matroids within the class of
transversal matroids.  We also contrast lattice path matroids with
fundamental transversal matroids and bicircular matroids.
Section~\ref{sec:high} treats higher connectivity.  The final section
introduces another minor-closed, dual-closed class of lattice path
matroids, the notch matroids, and characterizes this class by excluded
minors.

We assume familiarity with basic matroid theory (see,
e.g.,~\cite{ox,welsh}).  We follow the notation and terminology
of~\cite{ox}, with the following additions.  A flat $X$ of a matroid
$M$ is \emph{connected} if the restriction $M|X$ is connected.
A flat $X$ is \emph{trivial} if $X$ is independent; otherwise $X$ is
\emph{nontrivial}.  The flats in a collection $\mathcal{F}$ of flats
are \emph{incomparable}, or \emph{mutually incomparable}, if no flat
in $\mathcal{F}$ contains another flat in $\mathcal{F}$.  The
\emph{nullity}, $|X|-r(X)$, of a set $X$ is denoted by $\eta(X)$.
Recall that a matroid $M$ of rank $r$ is a \emph{paving matroid} if
every flat of rank less than $r-1$ is trivial.

Most matroids in this paper are transversal matroids
(see~\cite{bru,ing,welsh}). Recall that for a transversal matroid $M$,
a \emph{presentation} of $M$ is a multiset
$\mathcal{A}=(D_1,D_2,\ldots,D_k)$ of subsets of the ground set $E(M)$
such that the bases of $M$ are the maximal partial transversals of
$\mathcal{A}$.  As is justified by the following lemma
(see~\cite{bru}), we always consider presentations of rank-$r$
transversal matroids by set systems of size $r$.

\begin{lemma}\label{lem:minsets}
  Let $\mathcal{A}=(D_1,D_2,\ldots,D_k)$ be a presentation of a
  rank-$r$ transversal matroid $M$.  If some basis of $M$ is a
  transversal of $(D_{i_1},D_{i_2},\ldots,D_{i_r})$, with
  $i_1<i_2<\cdots<i_r$, then $(D_{i_1},D_{i_2},\ldots,D_{i_r})$ is
  also a presentation of $M$.
\end{lemma}

We use $[n]$ to denote the interval $\{1,2,\ldots,n\}$ of integers, and,
similarly, $[i,j]$ to denote the interval $\{i,i+1,\ldots,j\}$ of
integers.

\section{Background}\label{sec:back}

This section starts by reviewing the definition and basic properties
of lattice path matroids from~\cite{lpm1}.  The notation established
in this section is used throughout the paper. Also included are the
basic results about matroid connectivity that we use later.

Unless otherwise stated, all lattice paths in this paper start at the
point $(0,0)$ and use steps $E=(1,0)$ and $N=(0,1)$, which are called
\emph{East} and \emph{North}, respectively.  Paths are usually
represented as words in the alphabet $\{E,N\}$.  We say that a lattice
path $P$ has a \emph{$NE$ corner at $h$} if step $h$ of $P$ is North
and step $h+1$ is East. An \emph{$EN$ corner at $k$} is defined
similarly.  A corner can also be specified by the coordinates of the
point where the North and East steps meet.

A \emph{lattice path matroid} is, up to isomorphism, a matroid of the
type $M[P,Q]$ that we now define.  Let $P$ and $Q$ be lattice paths
from $(0,0)$ to $(m,r)$ with $P$ never going above $Q$.  Let
$\mathcal{P}$ be the set of all lattice paths from $(0,0)$ to $(m,r)$
that go neither above $Q$ nor below $P$.  For $i$ with $1\leq i\leq
r$, let $N_i$ be the set
$$N_i:=\{ j\, :\, \text{step $j$ is the $i$-th North step of some path
  in $\mathcal{P}$} \}.$$
Thus, $N_1, N_2,\ldots,N_r$ is a sequence
of intervals in $[m+r]$, and both the left endpoints and the right
endpoints form strictly increasing sequences; the left and right
endpoints of $N_i$ correspond to the positions of the $i$-th North
steps in $Q$ and $P$, respectively. The matroid $M[P,Q]$ is the
transversal matroid on the ground set $[m+r]$ that has
$(N_1,N_2,\ldots,N_r)$ as a presentation.  We call
$(N_1,N_2,\ldots,N_r)$ the \emph{standard presentation} of $M[P,Q]$.
Note that $M[P,Q]$ has rank $r$ and nullity $m$.

\begin{figure}
\begin{center}
  \epsfxsize 4.5truein \epsffile{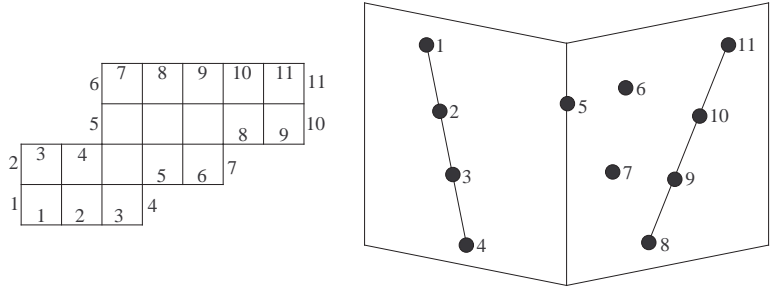}
\end{center}
\caption{A lattice path presentation and geometric representation of a
  lattice path matroid.}\label{draw}
\end{figure}

Figure~\ref{draw} shows a lattice path matroid of rank $4$ and nullity
$7$. The intervals in the standard presentation are $N_1=[4]$,
$N_2=[2,7]$, $N_3=[5,10]$, and $N_4=[6,11]$. (Section~\ref{ssec:geom}
explains how to find a geometric representation of a lattice path
matroid.)

A feature that enriches the subject of lattice path matroids is the
variety of ways in which these matroids can be viewed.  On the one
hand, the theory of transversal matroids provides many useful tools.
On the other hand, the following theorem from~\cite{lpm1} gives an
interpretation of the bases that leads to attractive descriptions of
many matroid concepts (see, e.g., \cite[Theorem 5.4]{lpm1} on basis
activities).

\begin{thm}\label{thm:bases}
  The map $R\mapsto \{i\,:\, \text{the } i\text{-th step of } R
  \text{ is North}\}$ is a bijection from $\mathcal{P}$ onto the set
  of bases of $M[P,Q]$.
\end{thm}

We use $\mathcal{L}$ to denote the class of all lattice path matroids.
We call the pair $(P,Q)$ a \emph{lattice path presentation} of
$M[P,Q]$, or, if there is no danger of confusion, a
\emph{presentation} of $M[P,Q]$.

Unless we say otherwise, all references to an order on the ground set
$[m+r]$ of $M[P,Q]$ are to the natural order $1<2<\cdots<m+r$.
However, this order is not inherent in the matroid structure; the
elements of a lattice path matroid typically can be linearly ordered
in many ways so as to correspond to the steps of lattice paths.  Also,
a lattice path matroid of rank $r$ and nullity $m$ need not have
$[m+r]$ as its ground set.  These comments motivate the following
definition.

\begin{dfn}\label{def:lpo}
  A linear ordering $s_1<s_2<\cdots<s_{m+r}$ of the ground set of a
  lattice path matroid $M$ is a \emph{lattice path ordering} if the
  map $s_i\mapsto i$ is an isomorphism of $M$ onto a lattice path
  matroid of the form $M[P,Q]$.
\end{dfn}

For some purposes it is useful to view lattice path matroids from the
following perspective, which does not refer to paths.  Lattice path
matroids are the transversal matroids $M$ for which $E(M)$ can be
linearly ordered so that $M$ has a presentation $(A_1,A_2,\ldots,A_r)$
where $A_i=[l_i,g_i]$ is an interval in $E(M)$ and the endpoints of these
intervals form two chains, $l_1<l_2<\cdots <l_r$ and $g_1<g_2<\cdots
<g_r$.

The \emph{incidence function} of a presentation $(A_1,A_2,\ldots,A_r)$
of a transversal matroid is given by $n(X) = \{i\,:\, X\cap
A_i\ne\emptyset\}$ for subsets $X$ of $E(M)$.  If no other
presentation is mentioned, the incidence function of the matroid
$M[P,Q]$ of rank $r$ and nullity $m$ is understood to be that
associated with the standard presentation.  For this incidence
function and for any element $x$ in $[m+r]$, the set $n(x)$ is an
interval in $[r]$; if $x<y$, then $\max\bigl(n(x)\bigr)\leq
\max\bigl(n(y)\bigr)$ and $\min\bigl(n(x)\bigr)\leq
\min\bigl(n(y)\bigr)$.

An independent set $I$ in a lattice path matroid $M[P,Q]$ is a partial
transversal of $(N_1, N_2,\ldots,N_r)$. Typically there are many ways
to match $I$ with $N_1, N_2,\ldots,N_r$.  The next two results show
that $I$ can always be matched in a natural way.  The following lemma,
which is crucial in the proof of Theorem~\ref{thm:bases}, is
from~\cite{lpm1}

\begin{lemma}\label{lem:basisnoncross}
  Assume $\{b_1,b_2,\ldots,b_r\}$ is a basis of a lattice path
  matroid $M[P,Q]$ with $b_1<b_2<\cdots<b_r$. Then $b_i$ is in $N_i$
  for all $i$.
\end{lemma}

Corollary~\ref{cor:indepnoncross} follows by extending the given
independent set $I$ to a basis and applying
Lemma~\ref{lem:basisnoncross}.

\begin{cor}\label{cor:indepnoncross}
  Assume $I$ is an independent set of a lattice path matroid $M[P,Q]$
  with $|I|=|n(I)|$; let $I$ be $\{a_1,a_2,\ldots,a_k\}$ with
  $a_1<a_2<\cdots<a_k$ and let $n(I)$ be $\{i_1,i_2,\ldots,i_k\}$ with
  $i_1<i_2<\cdots<i_k$. Then $a_j$ is in $N_{i_j}$ for all $j$ with
  $1\leq j \leq k$.
\end{cor} 

We now gather several results on matroid connectivity that are
relevant to parts of the paper.  The first result~\cite[Theorem
7.1.16]{ox} gives a fundamental link between connectivity and the
operation of parallel connection.

\begin{lemma}\label{lem:discon} 
  If $M$ is connected and $M/p$ is the direct sum $M_1\oplus M_2$,
  then $M$ is the parallel connection $P(M'_1,M'_2)$ of $M'_1:=M\del
  E(M_2)$ and $M'_2:=M\del E(M_1)$.
\end{lemma}

In Lemma~\ref{lem:discon}, since $M$ is connected, both $M'_1$ and
$M'_2$ are connected.  Recall that the rank
$r\bigl(P(M'_1,M'_2)\bigr)$ of a parallel connection whose basepoint
is not a loop is $r(M'_1)+r(M'_2)-1$.  These observations give the
following lemma.

\begin{lemma}\label{lem:smallerflats}
  If $M$ is connected, $x$ is not parallel to any element of $M$, and
  $M/x$ is disconnected, then there is a pair $A,B$ of nontrivial
  incomparable connected flats of $M$ with $r(A) + r(B)= r(M)+1$ and
  $A\cap B = \{x\}$.
\end{lemma}

The following useful lemma is easy to prove by using separating sets.

\begin{lemma}
  Assume that $X$ is a connected flat of a connected matroid $M$, that
  $x$ is in $X$, and that $M|(X-x)$ is connected.  Then $M\del x$ is
  connected.
\end{lemma}

The cyclic flats of a matroid $M$ (that is, the flats $F$ for which
$M|F$ has no isthmuses), together with their ranks, determine the
matroid~\cite[Proposition 2.1]{affine}.  As we show next, in the
loopless case it suffices to consider nontrivial connected flats.
Note that nontrivial connected flats are cyclic, but cyclic flats need
not be connected.  Thus, the next result is a mild refinement
of~\cite[Proposition 2.1]{affine}, and essentially the same idea
proves both results.

\begin{lemma}\label{lem:cyclicgen}
  The circuits of a loopless matroid $M$ (and hence $M$ itself) are
  determined by the nontrivial connected flats and their ranks.
\end{lemma}

\begin{proof}
  Note that if $C$ is an $i$-circuit, then $\cl(C)$ is a connected
  flat of rank $i-1$.  Thus, the circuits can be recovered inductively
  as follows: the $2$-circuits are the $2$-subsets of nontrivial
  rank-$1$ flats; the $3$-circuits are the $3$-subsets of $E(M)$ that
  contain no $2$-circuit and are subsets of connected lines, and so
  on.
\end{proof}

\begin{cor}\label{cor:genauto}
  The automorphisms of a loopless matroid are the permutations of the
  ground set that are rank-preserving bijections of the collection of
  nontrivial connected flats.
\end{cor}

\section{Basic Structural Properties of Lattice Path
  Matroids}\label{sec:lpm}

This section treats the basic structural properties of lattice path
matroids that play key roles throughout this paper.  Some of these
properties are shared by few other classes of matroids; for instance,
every nontrivial connected lattice path matroid has a spanning
circuit.  Other properties, such as the closure of the class of
lattice path matroids under minors and duals, while shared by many
classes of matroids, do not hold for the larger class of transversal
matroids.  Some of the properties are more technical and their
significance will become apparent only later in the paper. The topics
treated are fairly diverse, so we divide the material into subsections
that focus in the following issues: minors, duals, and direct sums;
connectivity and spanning circuits; the structure of circuits and
connected flats.

\subsection{Minors, Duals, and Direct Sums}\label{ssec:minor}

The class of transversal matroids, although closed under deletions and
direct sums, is closed under neither contractions nor duals.  In
contrast, we have the following result for lattice path matroids.

\begin{thm}\label{thm:minors}
  The class $\mathcal{L}$ is closed under minors, duals, and direct
  sums.
\end{thm}

  \begin{figure}
  \begin{center}
    \epsfxsize 2.75truein \epsffile{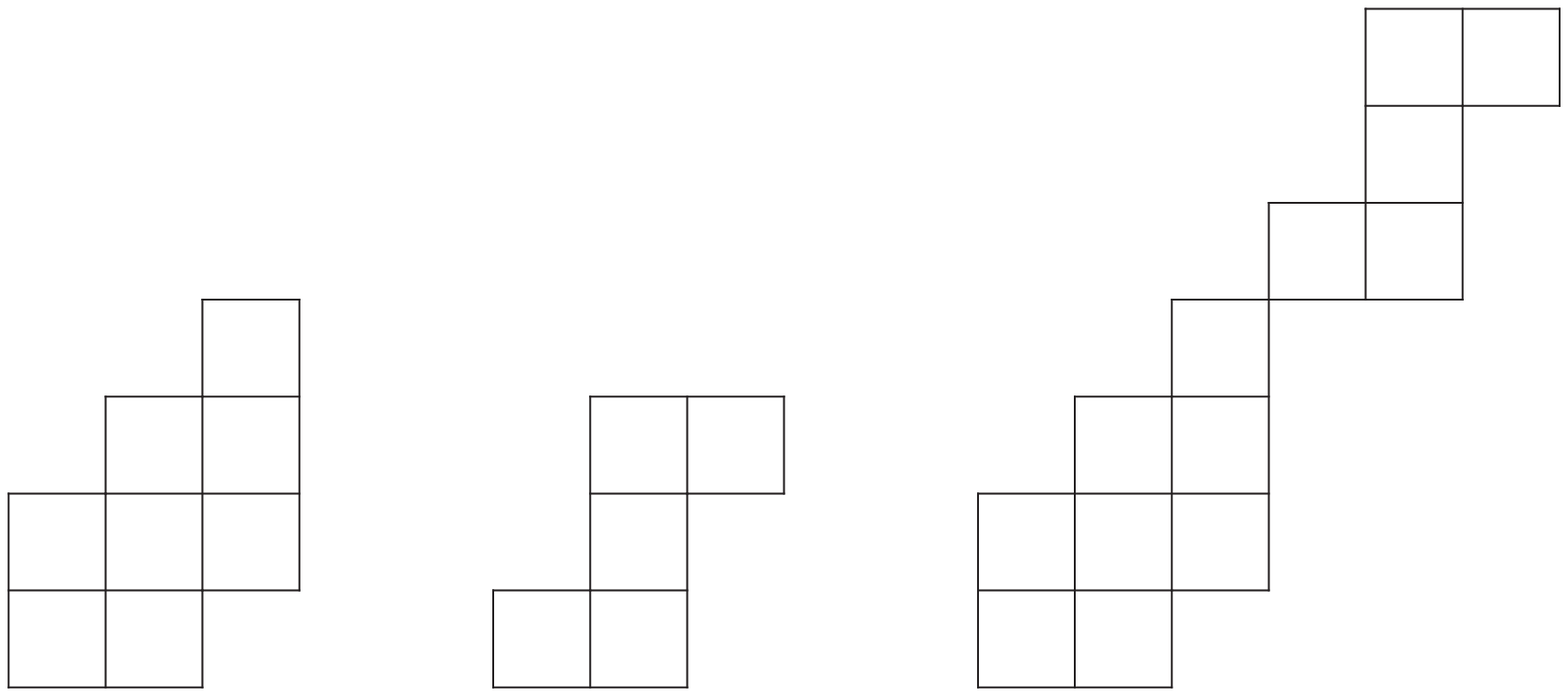}
  \end{center}
  \caption{Presentations of two lattice path matroids and their direct
    sum.}\label{sum}
  \end{figure}

\begin{proof}
  Figure~\ref{sum} illustrates the obvious construction to show that
  $\mathcal{L}$ is closed under direct sums.  For closure under
  duality, note that, from Theorem~\ref{thm:bases}, a basis of the
  dual of $M[P,Q]$ (i.e., the complement of a basis of $M[P,Q]$)
  corresponds to the East steps in a lattice path; the East steps of a
  lattice path are the North steps of the lattice path obtained by
  reflecting the entire diagram about the line $y=x$. This idea is
  illustrated in Figure~\ref{dual}.

  \begin{figure}
  \begin{center}
    \epsfxsize 2.7truein \epsffile{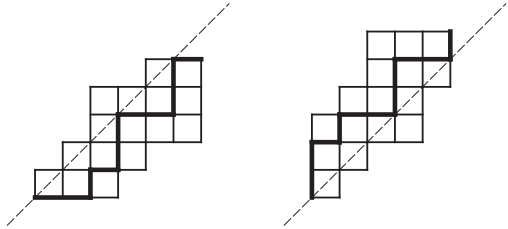}
  \end{center}
  \caption{Presentations of a lattice path matroid and its
    dual.}\label{dual}
  \end{figure}
  
  For closure under minors, it suffices to consider single-element
  deletions.  Let $x$ be in the lattice path matroid $M=M[P,Q]$ on
  $[m+r]$ with standard presentation $(N_1,N_2,\ldots,N_r)$.  Note
  that $(N_1-x,N_2-x,\ldots,N_r-x)$ is a presentation of $M\del x$;
  from this presentation, we will obtain one that shows that $M\del x$
  is a lattice path matroid.  Some set $N_i$ is $\{x\}$ if and only if
  $x$ is an isthmus of $M$; in this case, discard the empty set
  $N_i-x$ from the presentation above to obtain the required
  presentation of $M\del x$. Thus, assume $x$ is not an isthmus of
  $M$.  The sets $N_1-x,N_2-x,\ldots,N_r-x$ are intervals in the
  induced linear order on $[m+r]-x$.  In only two cases will the least
  elements or the greatest elements (or both) fail to increase
  strictly: (a) $x$ is the least element of the interval $N_i$ and
  $x+1$ is the least element of $N_{i+1}$, and (b) $x-1$ and $x$ are
  the greatest elements of $N_{j-1}$ and $N_j$, respectively. Assume
  case (a) applies.  Any basis of $M\del x$ (that is, any basis of $M$
  that does not contain $x$) that contains $x+1$ can, by
  Lemma~\ref{lem:basisnoncross}, be matched with $N_1,N_2,\ldots,N_r$
  so that $x+1$ is not matched to $N_{i+1}$.  Thus, the set system
  obtained by replacing $N_{i+1}$ by $N_{i+1}-\{x+1\}$ is also a
  presentation of $M\del x$. The same argument justifies replacing
  $N_{i+2}$ by $N_{i+2}-\{x+2\}$ if $x+2$ is the least element of
  $N_{i+2}$, and so on. Case (b) is handled similarly. The result is a
  presentation of $M\del x$ by intervals in which the least and
  greatest elements increase strictly, so $M\del x$ is a lattice path
  matroid.
\end{proof}

\begin{figure}
\begin{center}
  \epsfxsize 3.8truein \epsffile{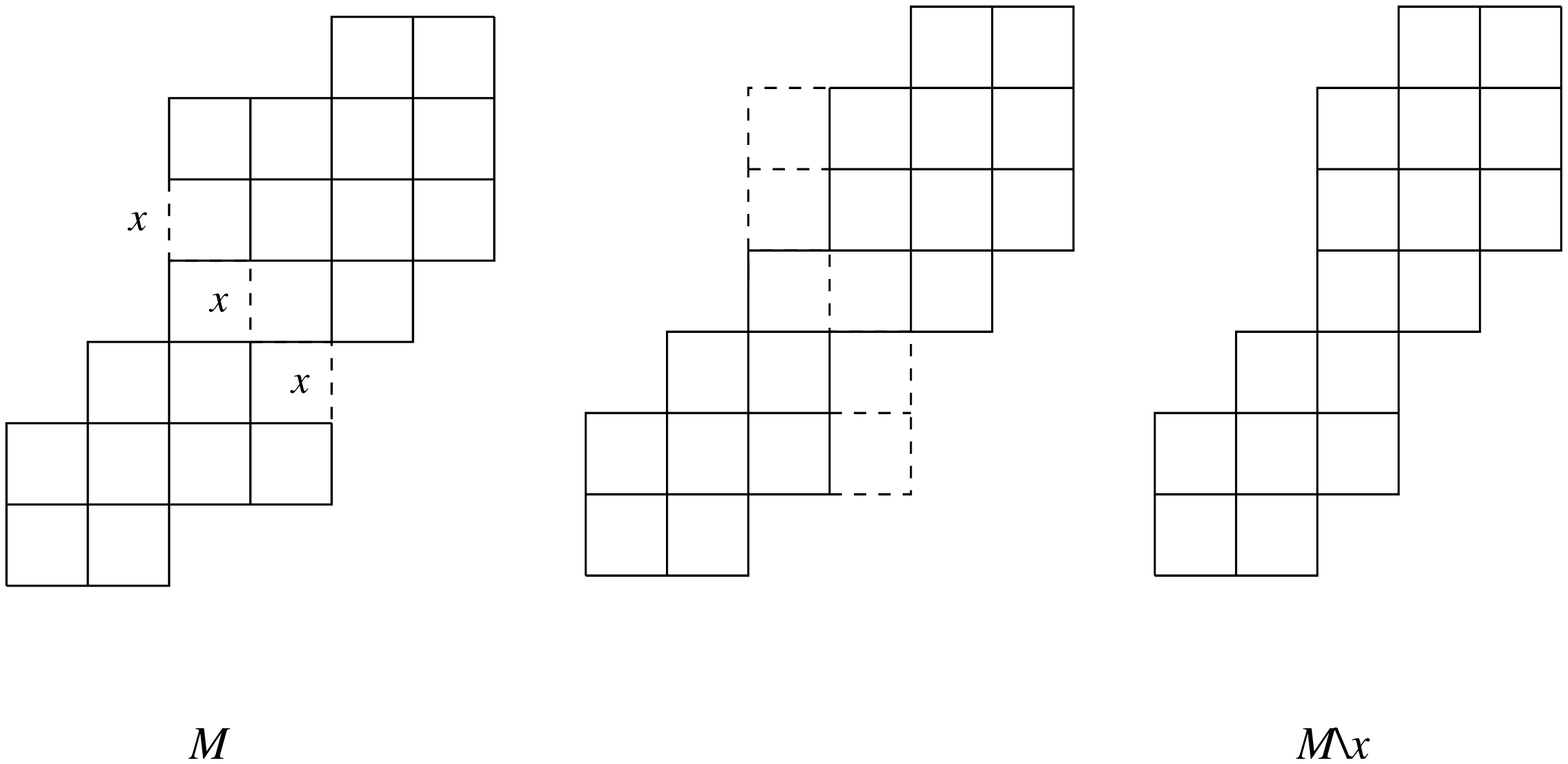}
\end{center}
\caption{The lattice path interpretation of the shortening of
  intervals that yields a presentation of a single-element
  deletion.}\label{delete}
\end{figure}

Single-element deletions and contractions can be described in terms of
the bounding paths of $M=M[P,Q]$ as follows. An isthmus is an element
$x$ for which some $N_i$ is $\{x\}$; to delete or contract $x$,
eliminate the corresponding common North step from both bounding
paths.  A loop is an element that is in no set $N_i$; to delete or
contract a loop, eliminate the corresponding common East step from $P$
and $Q$.  Now assume $x$ is neither a loop nor an isthmus. The upper
bounding path for $M\del x$ is formed by deleting from $Q$ the first
East step that is at or after step $x$; the lower bounding path for
$M\del x$ is formed by deleting from $P$ the last East step that is at
or before step $x$.  This is shown in Figure~\ref{delete}, where the
dashed steps in the middle diagram indicate the steps that bases of
$M\del x$ must avoid.  Dually, the upper bounding path for the
contraction $M/x$ is formed by deleting from $Q$ the last North step
that is at or before step $x$; the lower bounding path for $M/x$ is
formed by deleting from $P$ the first North step that is at or after
step $x$.

Corollary~\ref{cor:restint} treats restrictions of lattice path
matroids to intervals.  The lattice path interpretation of this result
is illustrated in Figure~\ref{fundint} on page~\pageref{fundint}.

\begin{cor}\label{cor:restint}
  Let $M$ be the lattice path matroid $M[P,Q]$ on the ground set
  $[m+r]$. Let $X$ be the initial segment $[i]$ and $Y$ be the final
  segment $[j+1,m+r]$ of $[m+r]$.  Let the $i$-th step of $Q$ end at
  the point $(h,k)$ and let the $j$-th step of $P$ end at $(h',k')$.
\begin{itemize}
\item[(a)] The bases of the restriction $M|X$ correspond to the
  lattice paths that go from $(0,0)$ to $(h,k)$ and go neither below
  $P$ nor above $Q$.
\item[(b)] The bases of the restriction $M|Y$ correspond to the
  lattice paths that go from $(h',k')$ to $(m,r)$ and go neither below
  $P$ nor above $Q$.
\item[(c)] If $h'\leq h$, then the bases of $M| (X\cap Y)$
  correspond to the lattice paths that go from $(h',k')$ to $(h,k)$
  and go neither below $P$ nor above $Q$.
\end{itemize}
\end{cor}

We close this section by noting that although $U_{1,2}\oplus
U_{1,2}\oplus U_{1,2}$ is a lattice path matroid, its truncation is
not transversal.  It follows that $\mathcal{L}$ is not closed under
the following operations: truncation, free extension, and elongation.

\subsection{Connectivity and Spanning Circuits}\label{ssec:conn}

We begin with a rare property.

\begin{thm}\label{thm:spanning}
  A connected lattice path matroid $M[P,Q]$ on at least two elements
  has a spanning circuit.
\end{thm}

\begin{proof}
  Let $M[P,Q]$ have rank $r$, let $N_j$ be $[l_j,g_j]$ for $1\leq j
  \leq r$, and let $C$ be the set
  $\{l_1,l_2,\ldots,l_{r-1},l_r,g_r\}$.  Showing that each set $C-x$,
  for $x$ in $C$, is a basis shows that $C$ is a spanning circuit.
  That $C-l_r$ and $C-g_r$ are bases is clear.  Since $M[P,Q]$ is not
  a direct sum of two matroids, $l_{i+1}$ must be in $N_i$ for $1\leq
  i <r$, from which it follows that each set $C-l_j$, with $1\leq j <
  r$, is a basis.
\end{proof}

It will be useful to single out the following immediate corollary of
Theorem~\ref{thm:spanning}.

\begin{cor}\label{cor:connflat}
  If $X$ is a nontrivial connected flat of a matroid $M$ and
  $M|X$ is a lattice path matroid, then $X$ is $\cl(C)$ for some
  circuit $C$ of $M$.
\end{cor}

The next theorem determines which lattice path matroids are connected.
One implication follows from the description of direct sums and the
other from the construction of the spanning circuit in the proof of
Theorem~\ref{thm:spanning}.

\begin{thm}\label{thm:conncrit}
  A lattice path matroid $M[P,Q]$ of rank $r$ and nullity $m$ is
  connected if and only if $P$ and $Q$ intersect only at $(0,0)$ and
  $(m,r)$.
\end{thm}

The parallel connection of two $3$-point lines, which has only one
spanning circuit, shows that there may be elements of a connected
lattice path matroid that are in no spanning circuit. There are
several ways to identify the elements of connected lattice path
matroids that are in spanning circuits.  The next result identifies
these elements via the standard presentation.

\begin{thm}\label{thm:spanningvar}
  An element $x$ of a nontrivial connected lattice path matroid
  $M[P,Q]$ of rank $r$ is in a spanning circuit of $M[P,Q]$ if and
  only if $x$ is in at least two of the sets $N_1,N_2,\ldots,N_r$, or
  $x$ is in $N_1$ or $N_r$.
\end{thm}

\begin{proof}
  Assume $x$ is in $N_i$ and $N_{i+1}$. Let $C$ be
  $\{l_1,l_2,\ldots,l_i,x,g_{i+1},g_{i+2},\ldots,g_r\}$ where $N_j$ is
  $[l_j,g_j]$.  By connectivity, we have $l_2\in N_1, l_3\in N_2,
  \ldots, l_i\in N_{i-1}$ and $g_{i+1}\in N_{i+2}, g_{i+2}\in N_{i+3},
  \ldots, g_{r-1}\in N_r$. An argument like that in the proof of
  Theorem~\ref{thm:spanning} shows that $C$ is a spanning circuit.
  Similar ideas show that $x$ is in a spanning circuit of $M[P,Q]$ if
  $x$ is in $N_1$ or $N_r$.
  
  Assume $n(x)$ is $\{i\}$ with $1<i<r$.  Note that the basepoint is
  in no spanning circuit of a parallel connection of matroids of rank
  two or more, so to complete the proof we need only show that
  $M[P,Q]$ is a parallel connection of two lattice path matroids, each
  of rank at least two, with basepoint $x$. Thus, by
  Lemma~\ref{lem:discon}, we need to show that $M[P,Q]/x\del X$ is
  disconnected where $X$ is the set of loops of $M[P,Q]/x$. This
  statement follows from the lattice path description of contraction
  along with the observations that $N_{i-1}$ contains only elements
  less than $x$ while $N_{i+1}$ contains only elements greater than
  $x$.
\end{proof}

The following characterizations of the elements that are in spanning
circuits use structural properties rather than presentations.

\begin{cor}\label{cor:spanningvar}
  Let $x$ be in a nontrivial connected lattice path matroid $M$.
\begin{itemize}
\item[(a)] No spanning circuit contains $x$ if and only if $M$ is a
  parallel connection of two lattice path matroids, each of rank at
  least two, with basepoint $x$.
\item[(b)] Some spanning circuit contains $x$ if and only if $M/x\del
  X$ is connected, where $X$ is the set of loops of $M/x$.
\end{itemize}
\end{cor}

\begin{proof}
  Part (a) follows from the proof of Theorem~\ref{thm:spanningvar}. If
  $x$ is in a spanning circuit $C$ of $M$, then $C-x$ is a spanning
  circuit of $M/x$, so $M/x\del X$ is connected. Conversely, if $x$ is
  in no spanning circuit of $M$, then, by part (a), $M$ is a parallel
  connection, with basepoint $x$, of matroids of rank at least two, so
  $M/x\del X$ is disconnected.
\end{proof}

\subsection{Circuits and Connected Flats}\label{ssec:flats}

Our first goal in this section is to characterize the circuits of
lattice path matroids.  This is done in Theorem~\ref{thm:cinterval},
the proof of which uses the following well-known elementary result
about the circuits of arbitrary transversal matroids.  This lemma
follows easily from Hall's theorem.

\begin{lemma}\label{lem:circform}
  Let $n$ be the incidence function of a presentation of a transversal
  matroid $M$.  If $C$ is a rank-$k$ circuit of $M$, then $|n(C)|$ is
  $k$, as is $|n(C-x)|$ for any $x$ in $C$.
\end{lemma}

\begin{thm}\label{thm:cinterval}
  Let $C = \{c_0,c_1,c_2,\ldots,c_k\}$ be a set in the lattice path
  matroid $M[P,Q]$; assume $c_0<c_1<c_2<\cdots<c_k$. Let $n(C)$ be
  $\{i_1,i_2,\ldots,i_s\}$, where $i_1<i_2<\cdots<i_s$. Then $C$ is a
  circuit of $M[P,Q]$ if and only if
\begin{itemize}
\item[($1$)] $s = k$,
\item[($2$)] $c_0\in N_{i_1}$,
\item[($3$)] $c_k\in N_{i_k}$, and
\item[($4$)] $c_j\in N_{i_j}\cap N_{i_{j+1}}$ for $j$ with $0<j<k$.
\end{itemize}
Furthermore, if $C$ is a circuit, then $i_{h+1} = i_h+1$ for $1\leq h
< k$.

\end{thm}

\begin{proof}
  It is immediate to check that if conditions ($1$)--($4$) hold, then
  $C$ is dependent and every $k$-subset of $C$ is a partial
  transversal and so is independent; thus $C$ is a circuit. For the
  converse, assume $C$ is a circuit. Assertion ($1$) follows from
  Lemma~\ref{lem:circform}, which also gives the equalities
  $|n(C-c_0)|=k=|n(C-c_k)|$.  Since $C-c_0$ is independent and
  $|n(C-c_0)|$ is $k$, it follows from
  Corollary~\ref{cor:indepnoncross} that $c_j$ is in $N_{i_j}$ for
  $1\leq j\leq k$. A similar argument using $C-c_k$ shows that $c_j$
  is in $N_{i_{j+1}}$ for $0\leq j\leq k-1$.  This proves assertions
  ($2$)--($4$). To prove the last assertion, assume there were an $h$
  not in $n(C)$ with $i_j < h < i_{j+1}$. From statement ($4$), we
  have that $c_j$ is in both $N_{i_j}$ and $N_{i_{j+1}}$. The
  inequalities
  $$\min(N_h)<\min (N_{i_{j+1}})\leq c_j\leq \max(N_{i_j})<\max(N_h)$$
  imply that $c_j$ is in $N_h$, which contradicts the assumption that
  $h$ is not in $n(C)$.
\end{proof}

By Lemma~\ref{lem:circform}, if $x$ is parallel to some element, then
$|n(x)|=1$.  By property ($4$) of Theorem~\ref{thm:cinterval}, at most
two elements $x$ in a circuit of a lattice path matroid can satisfy
the equality $|n(x)|=1$.  This observation proves the next result.

\begin{cor}\label{cor:cinterval}
  At most two elements in any circuit of a lattice path matroid are in
  nonsingleton parallel classes.
\end{cor}

The following result gives two useful properties of connected flats.

\begin{thm}\label{thm:finterval}
  Let $M[P,Q]$ have rank $r$ and nullity $m$.  Any nontrivial
  connected flat $X$ of $M[P,Q]$ is an interval in $[m+r]$ and $n(X)$
  is an interval of $r(X)$ elements in $[r]$.
\end{thm}

\begin{proof}
  The second assertion follows from Corollary~\ref{cor:connflat} and
  Theorem~\ref{thm:cinterval}.  For the first statement, let $n(X)$ be
  $[s,t]$ and assume $i<j<k$ with $i,k\in X$.  That $j$ is in $X$
  follows from the inequalities
  $$s \leq \min\bigl(n(i)\bigr) \leq \min\bigl(n(j)\bigr) \leq
  \max\bigl(n(j)\bigr) \leq \max\bigl(n(k)\bigr)\leq t.$$
\end{proof}

Theorem~\ref{thm:finterval} has many implications for the connected
flats of lattice path matroids, of which we mention four.

\begin{cor}\label{cor:finterval}
  Assume $M[P,Q]$ has rank $r$.
\begin{itemize}
\item[(i)] For $0\leq k\leq r-1$, there are at most $k+1$ nontrivial
  connected flats of rank $r-k$ in $M[P,Q]$.  In particular, $M[P,Q]$
  has at most two connected hyperplanes and at most $r-1$ connected
  lines.
\item[(ii)] A flat of positive rank of $M[P,Q]$ is covered by at
  most two connected flats.
\item[(iii)] The nontrivial connected flats of $M[P,Q]$ that are not
  contained in a fixed connected hyperplane $H$ of $M[P,Q]$ are
  linearly ordered by inclusion.
\item[(iv)] If $H$ and $H'$ are connected hyperplanes of $M[P,Q]$,
  then every nontrivial connected flat of $M[P,Q]$ is contained in at
  least one of $H$ and $H'$.
\end{itemize}
\end{cor}

The matroid $M[(E^2N)^{r-1}EN,NE(NE^2)^{r-1}]$, which is a parallel
connection of $r-1$ three-point lines in which elements have been
added parallel to the ``joints'' and the ``ends'', shows that all
upper bounds in parts (i) and (ii) of Corollary~\ref{cor:finterval}
are optimal.

The next result is another corollary of Theorem~\ref{thm:cinterval}.

\begin{cor}\label{cor:initfin}
  Let $C$ be the circuit $\{c_0,c_1,\ldots,c_k\}$ of $M[P,Q]$ with
  $c_0<c_1<\cdots<c_k$.  If $x$ is not in $C$ and $Z\cup x$ is a
  circuit of $M[P,Q]$ for some subset $Z$ of $C$, then $Z$ is either
  an initial segment $\{c_0,c_1,\ldots,c_i\}$ or a final segment
  $\{c_j,c_{j+1},\ldots,c_k\}$ of $C$.
\end{cor}

\begin{proof}
  The result follows from Lemma~\ref{lem:circform} and this simple
  corollary of Theorem~\ref{thm:cinterval}: for any proper subset $X$
  of $C$ that is neither an initial nor final segment of $C$, the
  inequality $|n(X)| > |X|$ holds.
\end{proof}

We conclude this section with a result we will use to show that
certain matroids are not lattice path matroids.

\begin{thm}\label{thm:notlpm}
  Assume a rank-$r$ matroid $M$ has two nontrivial connected flats $X$
  and $X'$ such that
\begin{itemize}
\item[($1$)] $X\cap X'\ne\emptyset$,
\item[($2$)] $r(X\cup X') = r$, and
\item[($3$)] $X\cup X'$ is a proper subset of the ground set $E(M)$ of $M$.
\end{itemize}
Then $M$ is not a lattice path matroid.
\end{thm}

\begin{proof}
  Assume, to the contrary, that $M$ is $M[P,Q]$. Fix $x$ in $X\cap X'$
  and $y$ in $E(M) -(X\cup X')$. By Theorem~\ref{thm:finterval}, along
  with assumptions ($1$) and ($2$), up to switching $X$ and $X'$ we
  would have $n(X)=[k]$ and $n(X')=[k',r]$ for some $k$ and $k'$ with
  $k'\leq k$. The inequality $y<x$ would give
  $\max\bigl(n(y)\bigr)\leq \max\bigl(n(x)\bigr)\leq k$, so $y$ would
  be in $\cl(X)$.  The inequality $x<y$ would give
  $\min\bigl(n(y)\bigr)\geq k'$, so $y$ would be in $\cl(X')$.  That
  these conclusions contradict the hypothesis proves the lemma.
\end{proof}

\section{Generalized Catalan Matroids}\label{sec:gcat}

Our next aim is to characterize lattice path matroids; this will be
done in Section~\ref{sec:fundamental}.  This section focuses on an
important subclass of $\mathcal{L}$ that has particularly simple
characterizations and many interesting properties.

\begin{dfn}
  The \emph{$n$-th Catalan matroid} $M_n$ is $M[E^nN^n,(EN)^n]$.  A
  \emph{generalized Catalan matroid} is, up to isomorphism, a matroid
  of the form $M[E^mN^r,Q]$.
\end{dfn}

For generalized Catalan matroids, the notation $M[P,Q]$ is simplified
to $M[Q]$.  We use $\mathcal{C}$ to denote the class of generalized
Catalan matroids. 

Generalized Catalan matroids have arisen in different contexts with a
corresponding variety of names and perspectives.  We gather here the
references currently known to us.  Crapo~\cite[Section 8]{crapofirst}
introduced these matroids to show that there are at least
$\binom{n}{r}$ nonisomorphic matroids of rank $r$ on $n$ elements.
His perspective was rediscovered in~\cite[Theorem 3.14]{lpm1}:
generalized Catalan matroids are precisely the matroids that are
obtained from the empty matroid by repeatedly applying the operations
of adding an isthmus and forming the free extension (this result is
generalized in Theorem~\ref{thm:extend} below).  By using ``nested''
presentations, Welsh~\cite{dw} proved that Crapo's lower bound on the
number of matroids holds within the smaller class of transversal
matroids.  These matroids arose again in~\cite{opr} in connection with
matroids defined in terms of integer-valued functions on finite sets.
They were studied further in~\cite{mixing}, where they were called
Schubert matroids and were shown to have the rapid mixing property.
In~\cite{ardila} they were rediscovered and related to shifted
complexes, and so acquired the name shifted matroids.  The link that
was established in~\cite{lpm1} between generalized Catalan matroids
and an enumerative problem known as the tennis ball problem influenced
the techniques used in~\cite{tbp} to solve that problem.
In~\cite{crapobill}, under the name of freedom matroids, generalized
Catalan matroids were used to construct a free algebra of matroids.

Catalan matroids have rich enumerative properties (see~\cite{lpm1}).
Their name comes from the fact that the number of bases of $M_n$ is
the Catalan number $C_n=\frac{1}{n+1}\binom{2n}{n}$; several other
invariants of $M_n$ are also Catalan numbers.  Although there is only
one Catalan matroid of each rank, these matroids generate the entire
class $\mathcal{C}$, in the sense of the following theorem.

\begin{thm}
  The smallest minor-closed class of matroids that contains all
  Catalan matroids is $\mathcal{C}$.
\end{thm}

\begin{proof}
  It follows from the lattice path interpretation of deletion and
  contraction given after the proof of Theorem~\ref{thm:minors} that
  $\mathcal{C}$ is closed under minors.  To see that any generalized
  Catalan matroid $M[Q]$ is a minor of a Catalan matroid, simply
  insert East and North steps into $Q$ so that the result is a Catalan
  matroid $M[(EN)^t]$. From $M[(EN)^t]$, delete the elements that
  correspond to the added East steps and contract the elements that
  correspond to the added North steps; by the lattice path
  interpretation of these operations, the resulting minor of
  $M[(EN)^t]$ is $M[Q]$.
\end{proof}

It is easy to see that $\mathcal{C}$, in addition to being closed
under minors, is closed under duals and (unlike $\mathcal{L}$) free
extension; therefore $\mathcal{C}$ is closed under truncation and
elongation. However, $\mathcal{C}$ is not closed under direct sums.

By Theorem~\ref{thm:conncrit}, a generalized Catalan matroid with at
least two elements is connected if and only if it has neither loops
nor isthmuses. The rest of this section focuses mainly on connected
generalized Catalan matroids since some results are slightly easier to
state with this restriction and, by what we just noted, there is
essentially no loss of generality.

The feature that makes generalized Catalan matroids easy to
characterize is the structure of the connected flats, as described in
the following lemma.

\begin{lemma}\label{lem:flatform}
  Assume $M[Q]$ has rank $r$, nullity $m$, and neither loops nor
  isthmuses.  Let the $EN$ corners of $Q$ be at steps
  $i_1,i_2,\ldots,i_k$ with $i_1<i_2<\cdots<i_k$.  The proper
  nontrivial connected flats of $M[Q]$ are the initial segments
  $[i_1]\subset[i_2]\subset\cdots\subset[i_k]$ of $[m+r]$.  The rank
  (resp.\ nullity) of $[i_h]$ is the number of North (resp.\ East)
  steps among the first $i_h$ steps of $Q$.
\end{lemma}

\begin{proof}
  The lemma follows easily once we show that any proper nontrivial
  connected flat $F$ of $M[Q]$ is an initial segment of $[m+r]$.  By
  Theorem~\ref{thm:finterval}, $F$ is an interval, say $[u,v]$, in
  $[m+r]$.  By Corollary~\ref{cor:restint}, the restriction of $M[Q]$
  to $[v]$ is $M[Q_v]$ where $Q_v$ consists of the first $v$ steps of
  $Q$.  Since $v$ is not an isthmus of $M[Q]|F$, it is not an isthmus
  of $M[Q_v]$, so the $v$-th step of $Q$ must be East.  Let $M[Q_v]$
  have rank $k$.  Note that $[v-k,v]$ is a spanning circuit of
  $M[Q_v]$ that is contained in $F$ and has closure $[v]$.  Thus, $F$
  is the initial segment $[v]$.
\end{proof}

The following result (which is essentially Lemma 2 of~\cite{opr}) is
an immediate corollary of Lemmas~\ref{lem:cyclicgen}
and~\ref{lem:flatform}.

\begin{cor}\label{cor:chainflats}
  A connected matroid is a generalized Catalan matroid if and only if
  its nontrivial connected flats are linearly ordered by inclusion.
\end{cor}

The following excluded-minor characterization of $\mathcal{C}$
from~\cite{opr} is not difficult to prove from
Corollary~\ref{cor:chainflats} and the results in
Section~\ref{sec:lpm}. Let $P_n$ be the truncation
$T_n(U_{n-1,n}\oplus U_{n-1,n})$ to rank $n$ of the direct sum of two
$n$-circuits.  Thus, $P_n$ is the paving matroid of rank $n$ whose
only nontrivial proper flats are two disjoint circuit-hyperplanes
whose union is the ground set.  It follows that $P_n$ is isomorphic to
$M[E^{n-1}NEN^{n-1},N^{n-1}ENE^{n-1}]$ and, by
Corollary~\ref{cor:chainflats}, that $P_n$ is not in $\mathcal{C}$.

\begin{thm}\label{thm:oxprow}
  A matroid is in $\mathcal{C}$ if and only if it has no minor
  isomorphic to $P_n$ for any $n\geq 2$.
\end{thm}

\section{Fundamental Flats and a Characterization of Lattice Path
  Matroids}\label{sec:fundamental}

While the structure of the connected flats of arbitrary connected
lattice path matroids is not as simple as that for generalized Catalan
matroids (Corollary~\ref{cor:chainflats}), this structure is still
easy to describe. We analyze this structure in this section and we use
it to characterize connected lattice path matroids.  We also show that
if $M[P,Q]$ is connected, then the paths $P$ and $Q$ are determined,
up to a $180^\circ$ rotation, by any matroid isomorphic to $M[P,Q]$.
The flats of central interest for these results are those we define
now.

\begin{dfn}
  Let $X$ be a connected flat of a connected matroid $M$ for which
  $|X|>1$ and $r(X)<r(M)$.  We say that $X$ is a \emph{fundamental
    flat of $M$} if for some spanning circuit $C$ of $M$ the
  intersection $X\cap C$ is a basis of $X$.
\end{dfn}

The first lemma shows how fundamental flats of lattice path matroids
reflect the order of the elements.

\begin{lemma}\label{lem:cflemma}
  Assume $M[P,Q]$ is connected and has rank $r$ and nullity $m$.  Let
  $X$ be a connected flat of $M[P,Q]$ with $|X|>1$ and $r(X)<r$. Then
  $X$ is a fundamental flat of $M[P,Q]$ if and only if $X$ is an
  initial or final segment of $[m+r]$.
\end{lemma}

\begin{proof}
  Let $N_i$ be $[l_i,g_i]$ for $1\leq i\leq r$.  If $X$ is an initial
  segment $[h]$ of $[m+r]$, then the spanning circuit
  $C=\{l_1,l_2,\ldots,l_r,g_r\},$ constructed in the proof of
  Theorem~\ref{thm:spanning}, has the property that $X\cap C$ is a
  basis of $X$.  Similarly, for a final segment $X$ of $[m+r]$, a
  spanning circuit with the required property is
  $\{l_1,g_1,g_2,\ldots,g_r\}$.
  
  Conversely, assume $C$ is a spanning circuit of $M[P,Q]$ and $X\cap
  C$ is a basis of $X$; say $C$ is $\{c_0,c_1,\ldots,c_r\}$ with
  $c_0<c_1<\cdots<c_r$. By Theorem~\ref{thm:finterval}, it suffices to
  show that either $1$ or $m+r$ is in $X$.  Let $x$ be in $X - C$ and
  let $C'$ be the unique circuit in $(X\cap C)\cup x$. By
  Corollary~\ref{cor:initfin}, $C'$ has the form
  $\{x,c_0,c_1,\ldots,c_u\}$ or $\{x,c_v,c_{v+1},\ldots,c_r\}$.  We
  will show that in the first case, $1$ is in $X$; a similar argument
  gives $m+r$ in $X$ in the second case. Thus, let $C'$ be
  $\{x,c_0,c_1,\ldots,c_u\}$. Note that $n(1)$ is $\{1\}$ and $1$ is
  in $n(c_0)$.  Note also that $\{c_0,c_1,\ldots,c_u\}$ is an
  independent set and, by Lemma~\ref{lem:circform} applied to $C'$, we
  have $|n(\{c_0,c_1,\ldots,c_u\})| = u+1$.  Thus,
  $$r(\{1,c_0,c_1,\ldots,c_u\}) \leq |n(\{1,c_0,c_1,\ldots,c_u\})|
  =u+1= r(\{c_0,c_1,\ldots,c_u\}).$$
  It follows that $1$ is in
  $\cl(\{c_0,c_1,\ldots,c_u\})$, so $1$ is in $X$, as claimed.
\end{proof}

Hence, to determine the fundamental flats of $M[P,Q]$, it suffices to
know which initial and final segments of $[m+r]$ are connected flats.
Note that the initial segment $[h]$ of $[m+r]$ is a proper nontrivial
connected flat, and hence a fundamental flat, if and only if the upper
path $Q$ has an $EN$ corner at $h$. Similarly, the final segment
$[k,m+r]$ of $[m+r]$ is a fundamental flat of $M[P,Q]$ if and only if
$P$ has a $NE$ corner at $k-1$.  These observations prove the
following theorem.

\begin{thm}\label{thm:shape}
  Assume $M[P,Q]$ is connected and has rank $r$ and nullity $m$.  Let
  the $EN$ corners of $Q$ be at $i_1,i_2,\ldots,i_h$, with
  $i_1<i_2<\cdots<i_h$, and the $NE$ corners of $P$ be at
  $j_1-1,j_2-1,\ldots,j_k-1$, with $j_1<j_2<\cdots<j_k$. The
  fundamental flats of $M[P,Q]$ are $[i_1]\subset
  [i_2]\subset\cdots\subset [i_h]$ and $[j_k,m+r]\subset \cdots
  \subset [j_2,m+r]\subset [j_1,m+r]$.
\end{thm}

Corollary~\ref{cor:twoch} follows immediately from
Theorem~\ref{thm:shape}.  Note that for generalized Catalan matroids,
the fundamental flats are precisely the flats given in
Lemma~\ref{lem:flatform}, so they form one chain under inclusion.

\begin{cor}\label{cor:twoch}
  The fundamental flats of a connected matroid $M$ in $\mathcal{L} -
  \mathcal{C}$ form two chains under inclusion; no set in one chain
  contains a set in the other chain.  Furthermore, for each pair $X,Y$
  of incomparable fundamental flats,
\begin{itemize}
\item[(a)] if $X\cap Y\ne\emptyset$, then $X\cup Y=E(M)$, and
\item[(b)] if $r(X)+r(Y)\geq r(M)$, then $r(X\cup Y) = r(M)$.
\end{itemize}
\end{cor}

While a connected lattice path matroid of rank $r$ has at most $k+1$
connected flats of rank $r-k$ (Corollary~\ref{cor:finterval}), it has
at most two fundamental flats of any given rank.

Theorem~\ref{thm:shape} and the lattice path interpretation of duality
give the next result.

\begin{cor}\label{cor:critcomp}
  For any lattice path matroid $M$, the fundamental flats of the dual
  $M^*$ are the set complements, $E(M)-F$, of the fundamental flats
  $F$ of $M$.
\end{cor}

A key observation that follows from Theorem~\ref{thm:shape} is that
although which flats are fundamental is independent of the order of
the elements that is inherent in any particular lattice path
presentation of a lattice path matroid, such a presentation makes it
easy to identify the fundamental flats.  Conversely, the chains of
fundamental flats give the bounding paths. More precisely, the paths
$P$ and $Q$ associated with $M[P,Q]$ are determined by the $NE$
corners of $P$ and the $EN$ corners of $Q$, and these corners are
determined by the ranks and nullities of the fundamental flats.
Typically there are two possible pairs of paths, according to which
chain of fundamental flats contains the least element of the ground
set.  These observations give the following theorem, which is one of
the main results of this section.

\begin{thm}\label{thm:shape2}
  The bounding paths $P$ and $Q$ of a connected lattice path matroid
  $M[P,Q]$ are determined by the matroid structure, up to a
  $180^\circ$ rotation.  That is, the only matroids $M[P^*,Q^*]$
  isomorphic to $M[P,Q]$ are $M[P,Q]$ and $M[Q^\rho,P^\rho]$ where
$(s_1s_2\cdots s_{m+r})^\rho$ is $s_{m+r}\cdots s_2 s_1$.
\end{thm}

Theorem~\ref{thm:shape} and its corollaries (including
Theorem~\ref{thm:shape2}) show that a connected lattice path matroid
is determined by its fundamental flats and their ranks.  The next
several results further develop this idea.  The following theorem
describes all connected flats of a connected lattice path matroid in
terms of its fundamental flats.

\begin{thm}\label{thm:cflats}
  Let $M$ be the connected lattice path matroid $M[P,Q]$ of rank $r$
  and nullity $m$ and let $F_1\subset F_2\subset\cdots\subset F_h$ and
  $G_1\subset G_2\subset\cdots\subset G_k$ be the chains of
  fundamental flats of $M$.  The proper nontrivial connected flats of
  $M$ are 
\begin{itemize}
\item[(i)] $F_1, F_2,\ldots, F_h, G_1, G_2,\ldots, G_k$, and
\item[(ii)] the intersections $F_i\cap G_j$ for which the inequality
  $m<\eta(F_i)+\eta(G_j)$ holds.
\end{itemize}
A nontrivial connected flat of the form $F_i\cap G_j$ has rank $r(F_i)
+ r(G_j) - r$.
\end{thm}

\begin{proof}
  The flats $F_1,F_2,\ldots,F_h,G_1,G_2,\ldots,G_k$, being
  fundamental, are connected.  The element $1$ is in either $F_1$ or
  $G_1$; we may assume it is in $F_1$.  For part (ii), we use
  Corollary~\ref{cor:restint} to find a lattice path presentation that
  shows that $F_i\cap G_j$ is connected. Using the notation in that
  corollary, let $X$ be $F_i$, so the point $(h,k)$ on $Q$ is
  $\bigl(\eta(F_i),r(F_i)\bigr)$; let $Y$ be $G_j$, so the point
  $(h',k')$ on $P$ is $\bigl(m-\eta(G_j),r-r(G_j)\bigr)$.  The
  inequality in part (ii) along with part (c) of
  Corollary~\ref{cor:restint} give a presentation of $M|(F_i\cap G_j)$
  (illustrated in Figure~\ref{fundint}) that, together with the fact
  that $P$ and $Q$ meet only at $(0,0)$ and $(m,r)$, implies that
  $F_i\cap G_j$ is connected and nontrivial.

  \begin{figure}
  \begin{center}
   \epsfxsize 3.6truein \epsffile{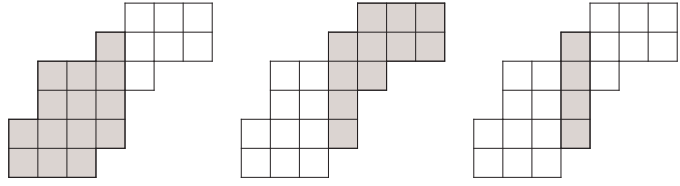}
  \end{center}
  \caption{The shaded regions show the 
    presentations of $M|F_i$, $M|G_j$, and $M|(F_i\cap
    G_j)$.}\label{fundint}
  \end{figure}
  
  Now assume $X$ is a proper nontrivial connected flat.  By
  Theorem~\ref{thm:finterval}, $X$ is an interval, say $[u,v]$, in
  $[m+r]$.  As in the proof of Lemma~\ref{lem:flatform}, it follows
  that the $u$-th step of $P$ and the $v$-th step of $Q$ are East
  steps.  Since $X$ is a flat, both $r(X\cup\{u-1\})$ and
  $r(X\cup\{v+1\})$ exceed $r(X)$, so step $u-1$ of $P$ and step $v+1$
  of $Q$, if there are such steps, are North steps.  From these
  observations and Theorem~\ref{thm:shape}, it follows that $X$ is of
  the form $F_i$, $G_j$, or $F_i\cap G_j$.  We need to show that if
  $F_i\cap G_j$ is connected, then the inequality
  $m-\eta(G_j)<\eta(F_i)$ holds.  This inequality follows by viewing
  $M|(F_i\cap G_j)$ as a restriction of $M|F_i$ and using the path
  presentations of these matroids given in
  Corollary~\ref{cor:restint}.  Indeed, from the lattice path diagrams
  (Figure~\ref{fundint}) it follows that $M|(F_i\cap G_j)$ is either
  free or connected, and the latter holds precisely when
  $\bigl(m-\eta(G_j),r-r(G_j)\bigr)$ is strictly to the left of
  $\bigl(\eta(F_i),r(F_i)\bigr)$.

  Lastly, let the connected flat $X$ be $F_i\cap G_j$.  From lattice
  path diagrams, we get
  $$r(M) = \bigl(r(F_i) - r(X)\bigr) + \bigl(r(G_j) - r(X)\bigr) +
  r(X),$$
  from which the last assertion follows.
\end{proof}

It follows from Theorem~\ref{thm:cflats} that any intersection of
connected flats is either a fundamental flat or an intersection of two
fundamental flats.  From this observation and the second paragraph of
the proof, it follows that a nonempty intersection of connected flats
is either connected or trivial.  Despite what the last part of
Theorem~\ref{thm:cflats} might suggest, it is easy to construct
examples in which the fundamental flats of lattice path matroids are
not modular.

The image, under an automorphism, of a fundamental flat of any matroid
is also fundamental.  This observation, Corollary~\ref{cor:genauto},
and Theorem~\ref{thm:cflats} give the following result.

\begin{cor}\label{cor:lpmauto}
  The automorphisms of a connected lattice path matroid are the
  permutations of the ground set that are rank-preserving bijections
  of the collection of fundamental flats.
\end{cor}

The proof of the second main result of this section,
Theorem~\ref{thm:lpmchar}, uses the following basic notions about
ordered sets. A \emph{strict partial order} is an irreflexive,
transitive relation.  Thus, strict partial orders differ from partial
orders only in whether each element is required to be unrelated, or
required to be related, to itself.  Given a strict partial order $<$
on $S$, elements $x$ and $y$ of $S$ are \emph{incomparable} if neither
$x<y$ nor $y<x$ holds.  \emph{Weak orders} are strict partial orders
in which incomparability is an equivalence relation.  Thus, linear
orders are weak orders in which the incomparability classes are
singletons.  Two weak orders $<_1$ and $<_2$ on $S$ are
\emph{compatible} if whenever elements $x$ and $y$ of $S$ are
comparable in both $<_1$ and $<_2$, and $x<_1 y$, then $x<_2 y$.

\begin{lemma}\label{lem:weak}
  Any two compatible weak orders have a common linear extension.
\end{lemma}

\begin{proof}
  Let $<_1$ and $<_2$ be compatible weak orders on $S$ and let the
  relation $<$ on $S$ be defined as follows: $x < y$ if either $x<_1
  y$ or $x<_2 y$. It is easy to check that $<$ is a weak order.  The
  lemma follows since $<$, like any strict partial order, can be
  extended to a linear order.
\end{proof}

We now turn to the second main result of the section.  This theorem
shows that the properties we developed above for the fundamental flats
and the connected flats of connected lattice path matroids
characterize these matroids.

\begin{thm}\label{thm:lpmchar}
  A connected matroid $M$ is a lattice path matroid if and only if the
  following properties hold.
\begin{itemize}
\item[(i)] The fundamental flats form at most two disjoint chains
  under inclusion, say $F_1\subset F_2\subset\cdots\subset F_h$ and
  $G_1\subset G_2\subset\cdots\subset G_k$.
\item[(ii)] If $F_i\cap G_j\ne \emptyset$, then $F_i\cup G_j = E(M)$.
\item[(iii)] The proper nontrivial connected flats of $M$ are
  precisely the following sets:
\begin{itemize}
\item[(a)] $F_1, F_2,\ldots, F_h, G_1, G_2,\ldots, G_k$, and
\item[(b)] intersections $F_i\cap G_j$ for which the inequality
  $m<\eta(F_i)+\eta(G_j)$ holds.
\end{itemize}
\item[(iv)] The rank of the flat $F_i\cap G_j$ of item \emph{(iii:b)}
  is $r(F_i)+r(G_j)-r(M)$
\end{itemize}
\end{thm}

\begin{proof}
  By Theorem~\ref{thm:shape}, Lemma~\ref{lem:flatform}, and
  Corollary~\ref{cor:chainflats}, $M$ is a generalized Catalan matroid
  if and only if properties (i)--(iv) hold where there is at most one
  chain of fundamental flats. By Theorems~\ref{thm:shape}
  and~\ref{thm:cflats}, the fundamental flats of a lattice path
  matroid that is not a generalized Catalan matroid satisfy properties
  (i)--(iv) with neither chain of fundamental flats being empty. Hence
  we need only prove the converse in the case that neither chain of
  fundamental flats is empty.
  
  Assume $M$ has rank $r$ and nullity $m$.  To show that $M$ is a
  lattice path matroid, we construct lattice paths $P$ and $Q$ and an
  isomorphism of $M$ onto $M[P,Q]$.  To show that $P$ stays strictly
  below $Q$ except at $(0,0)$ and $(m,r)$, we will use the following
  statements about fundamental flats.
  \begin{itemize}
    \item[(A)] If $F_i\cap G_j\ne\emptyset$, then $r(F_i)+r(G_j) >r$.
    \item[(B)] If $F_i\cap G_j=\emptyset$, then $\eta(F_i)+\eta(G_j) <m$.
  \end{itemize}
  To prove statement (A), note that we have the inequality
  $$
  r(F_i)+r(G_j) \geq r(F_i\cup G_j)+r(F_i\cap G_j) = r(M)+
  r(F_i\cap G_j)$$
  by semimodularity and property (ii).  Since $M$ has
  no loops, $r(F_i\cap G_j)$ is positive, so the desired inequality
  follows.  To prove statement (B), first recall that $\eta$ is
  nondecreasing, i.e., if $X\subseteq Y$, then $\eta(X)\leq \eta(Y)$.
  Since $F_i$ and $G_j$ are disjoint, we have $\eta(F_i) + \eta(G_j) =
  |F_i \cup G_j|-r(F_i) - r(G_j)$.  Thus, if $r(F_i) +
  r(G_j)>r(F_i\cup G_j)$, then we have $\eta(F_i) + \eta(G_j) <
  \eta(F_i\cup G_j)\leq m$. If $r(F_i) + r(G_j)=r(F_i\cup G_j)$, then
  $M|(F_i\cup G_j)$ is disconnected and we have the equality
  $\eta(F_i) + \eta(G_j)=\eta(F_i \cup G_j)$.  Since $M$ is connected,
  we have $\eta(F_i \cup G_j)<\eta(M)$, which gives the desired
  inequality.
  
  Let lattice paths $P$ and $Q$ from $(0,0)$ to $(m,r)$ be given as
  follows.
  \begin{itemize}
  \item[(a)] The $NE$ corners of $P$ are at the points
    $\bigl(m-\eta(G_j), r-r(G_j)\bigr)$ for $j$ in $[k]$.
  \item[(b)] The $EN$ corners of $Q$ are at the points
    $\bigl(\eta(F_i), r(F_i)\bigr)$ for $i$ in $[h]$.
   \end{itemize}
   Note that $P$ stays strictly below $Q$ except at the endpoints if
   and only if for every $NE$ corner $(x_P,y_P)$ of $P$ and every $EN$
   corner $(x_Q,y_Q)$ of $Q$, at least one of the inequalities
   $x_Q<x_P$ and $y_Q>y_P$ holds.  These inequalities are those in
   statements (A) and (B), so $P$ stays strictly below $Q$ except at
   $(0,0)$ and $(m,r)$.
   
   To construct an isomorphism of $M$ onto $M[P,Q]$, we define a
   linear order on $E(M)$ that we use to map $E(M)$ onto $[m+r]$, the
   ground set of $M[P,Q]$. We first define two relations $<_F$ and
   $<_G$ on $E(M)$.  Let $F_{h+1}$ and $G_{k+1}$ be $E(M)$. Define
   $<_F$ as follows: $x<_F y$ for $x,y\in E(M)$ if there is an integer
   $i$ in $[h]$ with $x\in F_i$ and $y\in F_{i+1}-F_i$.  Note that
   $<_F$ is a weak order whose incomparability classes are $F_1$ and
   the set differences $F_{i+1}-F_i$. Define $<_G$ similarly: $x<_G y$
   for $x,y\in E(M)$ if there is an integer $j$ in $[k]$ with $x\in
   G_{j+1}-G_j$ and $y\in G_j$. Thus, $<_G$ is also a weak order and
   the incomparability classes are $G_1$ and the differences
   $G_{i+1}-G_i$.  Note that if we had $x <_F y$ and $y <_G x$, then
   there would be fundamental flats $F_i$ and $G_j$ that both contain
   $x$ and not $y$, contrary to hypothesis (ii).  Thus, the weak
   orders $<_F$ and $<_G$ are compatible, so by Lemma~\ref{lem:weak}
   there is a linear order, say $x_1<x_2<\cdots<x_{m+r}$, of $E(M)$
   that extends both $<_F$ and $<_G$.
   
   Let $\phi:E(M) \rightarrow [m+r]$ be given by $\phi(x_i) = i$.  By
   construction, $\phi$ is a bijection of $E(M)$ onto $[m+r]$ that is
   a rank-preserving bijection of the fundamental flats of $M$ onto
   the fundamental flats of $M[P,Q]$. Furthermore, by assumptions
   (iii) and (iv) and Theorem~\ref{thm:cflats}, $\phi$ is a
   rank-preserving bijection of the set of connected flats of $M$ onto
   those of $M[P,Q]$.  By Lemma~\ref{lem:cyclicgen}, it follows that
   $\phi$ is an isomorphism of $M$ onto $M[P,Q]$; thus, $M$ is a
   lattice path matroid.
\end{proof}

We close this section by giving a pair of six-element matroids that
have the same collection of fundamental flats, yet only one of which
is in $\mathcal{L}$; thus, conditions (i) and (ii) in
Theorem~\ref{thm:lpmchar} are not enough to characterize lattice path
matroids.  The uniform matroid $U_{4,6}$ is a lattice path matroid
with no fundamental flats since the bounding paths are $P=E^2N^4$ and
$Q=N^4E^2$.  The prism (the matroid $C_{4,2}$ of Figure~\ref{b22} on
page~\pageref{b22}) is not a lattice path matroid (condition (iii) of
Theorem~\ref{thm:lpmchar} fails) and, since it has no spanning
circuits, it too has no fundamental flats.

\section{Lattice Path Matroids as Transversal Matroids}\label{sec:trans}

The aspects of lattice path matroids treated in this section relate to
important topics in the theory of transversal matroids.  We start by
characterizing the set systems that are maximal presentations of
lattice path matroids.  This result plays a key role in an algorithm
for determining whether a transversal matroid is in $\mathcal{L}$.  By
combining the result on maximal presentations with Brylawski's affine
representation of transversal matroids, we get a geometric description
of lattice path matroids.  We conclude the section by comparing
$\mathcal{L}$ with the dual-closed class of fundamental transversal
matroids and the minor-closed class of bicircular matroids.

\subsection{Maximal and Minimal Presentations}\label{ssec:pres}

Two types of presentations are of interest in this section.  A
presentation $\mathcal{A}=(A_1,A_2,\ldots,A_r)$ of a transversal
matroid $M$ is \emph{minimal} if the only presentation
$(A'_1,A'_2,\ldots,A'_r)$ of $M$ with $A'_i$ contained in $A_i$ for
all $i$ is $\mathcal{A}$.  The presentation $\mathcal{A}$ is
\emph{maximal} if the only presentation $(A'_1,A'_2,\ldots,A'_r)$ of
$M$ with $A_i$ contained in $A'_i$ for all $i$ is $\mathcal{A}$.  It
is well known that while each transversal matroid has a unique maximal
presentation, it typically has many minimal presentations.
(See, e.g.,~\cite{bondy,bru,ing}.)

\begin{thm}
  Standard presentations of lattice path matroids are minimal.
\end{thm}

\begin{proof}
  Let $(N_1,N_2,\ldots,N_r)$ be the standard presentation of the
  matroid $M[P,Q]$ and let $(N'_1,N'_2,\ldots,N'_r)$ be any
  presentation of $M[P,Q]$ with $N'_i\subseteq N_i$ for $1\leq i\leq
  r$.  To prove the theorem, we must show the inclusion $N_i\subseteq
  N'_i$ for $1\leq i\leq r$.  Let $x$ be in $N_i$.  Let $B$ consist of
  the least elements of $N_1,N_2,\ldots,N_{i-1}$, the greatest
  elements of $N_{i+1},N_{i+2},\ldots,N_{r}$, and $x$.  Thus, $B$ is a
  basis of $M[P,Q]$. Note that for $B$ to be a transversal of
  $(N'_1,N'_2,\ldots,N'_r)$, the element $x$ must be in $N'_i$, as
  needed.
\end{proof}

With the following result of Bondy~\cite{bondy}, we will get a simple
description, in terms of intervals, of the maximal presentation of a
lattice path matroid.

\begin{lemma}\label{lem:max}
  Given a presentation $(A_1,A_2,\ldots,A_r)$ of a rank-$r$
  transversal matroid $M$, the maximal presentation of $M$ is
  $(A_1\cup I_1,A_2\cup I_2,\ldots,A_r\cup I_r)$ where $I_j$ is the
  set of isthmuses of the deletion $M\del A_j$.
\end{lemma}

Together with Lemma~\ref{lem:max}, the following result from~\cite{hk}
implies that from any presentation of a transversal matroid, the
maximal presentation can be found in polynomial time in the size of
the ground set.  This observation will be important in the algorithm
for recognizing lattice path matroids among transversal matroids.

\begin{lemma}\label{lem:poly}
  The maximal size of a matching in a bipartite graph can be found in
  polynomial time in the number of vertices.
\end{lemma}

The discussion below focuses on matroids that have no isthmuses. This
restriction is justified by noting that the isthmuses of a transversal
matroid are in all sets in the maximal presentation, and so are easy
to deal with.

Let $(N_1,N_2,\ldots,N_r)$ be the standard presentation of the lattice
path matroid $M=M[P,Q]$ on $[m+r]$, where $M$ has no isthmuses. Let
$N_i$ be $[l_i,g_i]$.  Theorem~\ref{thm:cinterval} implies that each
connected component of $M\del N_i$ is a subset of either $[g_i+1,m+r]$
or $[l_i-1]$.  Thus, the set of isthmuses of $M\del N_i$ is the union
of the sets $I_i^+$ and $I_i^-$ of isthmuses of the restrictions of
$M$ to $[g_i+1,m+r]$ and $[l_i-1]$, respectively.
Corollary~\ref{cor:restint} implies that $I_i^+$ and $I_i^-$ are given
as follows:
\begin{equation}\label{eq:isth1}
I_i^+=\{g_i+j \,:\, g_i+j \text{ is the greatest element of }
N_{i+j}, j>0\},
\end{equation}
\begin{equation}\label{eq:isth2}
I_i^-=\{l_i-j \,:\, l_i-j \text{ is the least element of }
N_{i-j}, j>0\}.
\end{equation}
This proves the following theorem.

\begin{figure}
\begin{center}
  \epsfxsize 4.0truein \epsffile{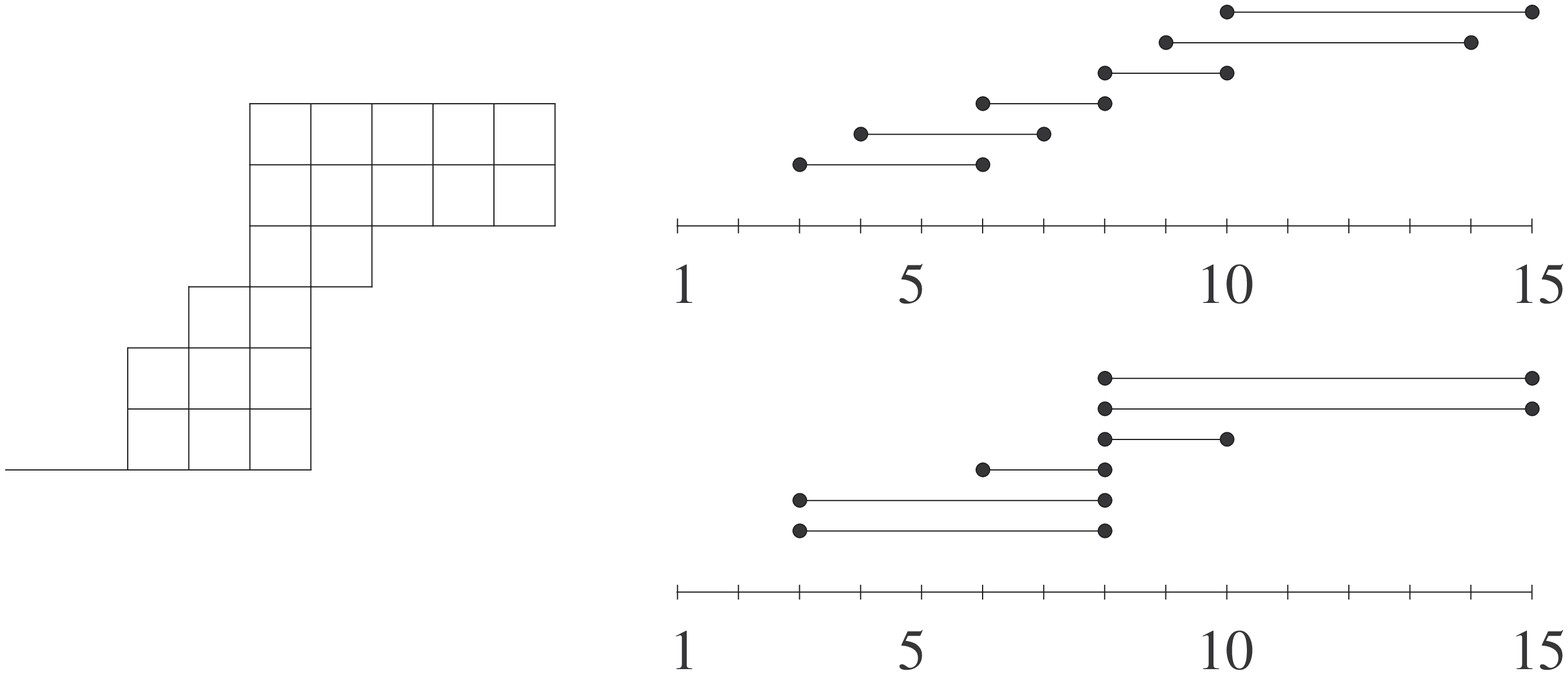}
\end{center}
\caption{The standard and maximal presentations of a lattice path
  matroid.}\label{maximal} 
\end{figure}

\begin{thm}\label{thm:max}
  Let $(N_1,N_2,\ldots,N_r)$ be the standard presentation of the
  lattice path matroid $M[P,Q]$ that has no isthmuses.  The maximal
  presentation of $M[P,Q]$ is $(N'_1,N'_2,\ldots,N'_r)$ where $N'_i$
  is $N_i\cup I_i^+ \cup I_i^-$ and $I_i^+$ and $I_i^-$ are given by
  Eqs.~(\ref{eq:isth1})--(\ref{eq:isth2}).
\end{thm}

The sets in the maximal presentation of a lattice path matroid have a
simple graphical interpretation, as Figure~\ref{maximal} illustrates.
While there are no containments among intervals in the standard
presentation, this figure shows that there may be containments (even
equalities) among intervals in the maximal presentation.

Theorem~\ref{thm:charint}, which characterizes the multisets of
intervals in $[m+r]$ that are maximal presentations of lattice path
matroids, uses the following notation.  For an indexed multiset
$(T_1,T_2,\ldots,T_r)$ of nonempty intervals in $[m+r]$ with
$T_i=[a_i,b_i]$, write $T_i\prec T_j$ if either $a_i<a_j$ or
$b_i<b_j$.  Thus, two intervals are unrelated if and only if they are
equal.  For arbitrary multisets of intervals, both $T_i\prec T_j$ and
$T_j\prec T_i$ may hold; in contrast, if $(T_1,T_2,\ldots,T_r)$ is the
maximal presentation of a lattice path matroid, then $\prec$ is a weak
order.  If $\prec$ is a weak order, then we assume that the set system
$(T_1,T_2,\ldots,T_r)$ is indexed so that we can have $T_i\prec T_j$
only for $i<j$.  In this case, let $d(T_h)$ be
$|\{i\,:\,i<h,a_i=a_h\}|$ and let $d'(T_h)$ be
$|\{j\,:\,h<j,b_h=b_j\}|$.

\begin{thm}\label{thm:charint}
  A set system $(T_1,T_2,\ldots,T_r)$ of nonempty intervals in $[m+r]$
  is the maximal presentation of a rank-$r$ lattice path matroid on
  $[m+r]$ that has no isthmuses if and only if 
\begin{itemize}
\item[(i)] the relation $\prec$ is a weak order,
\item[(ii)] for all pairs $T_i$ and $T_j$, neither $|T_i-T_j|$ nor
  $|T_j-T_i|$ is $1$, and
\item[(iii)] $d(T_i) + d'(T_i) + 2\leq |T_i|$ for every $i$.
\end{itemize}
\end{thm}

\begin{proof}
  For the maximal presentation of a lattice path matroid $M[P,Q]$ with
  no isthmuses, properties (i)--(iii) follow from
  Theorem~\ref{thm:max}.  For the converse, note that removing from
  $T_i$ its least $d(T_i)$ elements and its greatest $d'(T_i)$ yields
  the standard presentation of a lattice path matroid that, by
  property (iii), has no isthmuses and for which
  $(T_1,T_2,\ldots,T_r)$ is, by Theorem~\ref{thm:max}, the maximal
  presentation.
\end{proof}

\subsection{Recognizing Lattice Path Matroids}\label{ssec:recog}

When treating algorithmic questions about matroids, it is usual to
assume that a matroid is given by an \emph{independence oracle}, that
is, a subroutine that outputs, in constant time, whether a subset of
the ground set is independent.  While there are algorithms that
recognize transversal matroids within the class of all matroids
(see~\cite{bd}), Jensen and Korte~\cite{jensen} have shown that there
is no polynomial-time algorithm to decide if a matroid is transversal
from an independence oracle. The same proof as in~\cite{jensen} shows
that there is no such algorithm to decide whether a matroid is a
lattice path matroid.  Transversal matroids are more conveniently
specified by set systems than by independence oracles.  This section
gives a polynomial-time algorithm that, given a set system, decides
whether the corresponding transversal matroid is a lattice path
matroid.

We start with some simplifications.  A presentation $\mathcal{A}$ of
$M$ can be represented by a bipartite graph $\Delta[\mathcal{A}]$ in
the obvious way \cite[Section 1.6]{ox}.  Therefore, by
Lemma~\ref{lem:poly}, the isthmuses of a transversal matroid can be
identified and deleted in polynomial time.  If $M$ has no isthmuses,
then the connected components of $M$ come from those of
$\Delta[\mathcal{A}]$.  These observations and
Theorem~\ref{thm:minors} justify focusing on connected transversal
matroids.  As noted in Section~\ref{ssec:pres}, the maximal
presentation can be found from any presentation in polynomial time, so
we focus on maximal presentations.

The key to the recognition algorithm below is to efficiently recover
lattice path orderings from the maximal presentation.  We begin with
some observations that relate these notions.  Assume
$\mathcal{A}=(A_1,A_2,\ldots,A_r)$ is the maximal presentation of the
connected lattice path matroid $M[P,Q]$ on the ground set $[m+r]$ and
let $n$ be the incidence function of $\mathcal{A}$.  Let
$C_1,C_2,\ldots,C_k$ be the equivalence classes of the relation on
$[m+r]$ in which $x$ and $y$ are related if and only if $n(x)=n(y)$.
Each set $C_i$ is an interval in $[m+r]$. We may assume that
$C_1,C_2,\ldots,C_k$ are indexed so that $x_1<x_2<\cdots<x_k$ for any
elements $x_1,x_2,\ldots,x_k$ with $x_i$ in $C_i$.  Since $M[P,Q]$ is
connected, we have $n(C_i)\cap n(C_{i+1})\ne\emptyset$ for $i$ with
$1\leq i <k$. Any permutation $\sigma$ of $[m+r]$ with
$\sigma(C_i)=C_i$ for $1\leq i\leq k$ is clearly an automorphism of
$M[P,Q]$, so the linear order $\sigma(1)<\sigma(2)<\cdots<\sigma(m+r)$
is a lattice path order, as is
$\sigma(m+r)<\cdots<\sigma(2)<\sigma(1)$. Relative to any of these
linear orders, the sets in $\mathcal{A}$ are intervals and the
properties in Theorem~\ref{thm:charint} hold. These lattice path
orderings of $[m+r]$ are essentially equivalent to the orderings
$C_1<C_2<\cdots<C_k$ and $C_k<C_{k-1}<\cdots<C_1$ of
$C_1,C_2,\ldots,C_k$. Observe that $C_1,C_2,\ldots,C_k$ and
$C_k,C_{k-1},\ldots,C_1$ are the only permutations
$X_1,X_2,\ldots,X_k$ of $C_1,C_1,\ldots,C_k$ that satisfy the
following property.
\begin{quote} \textbf{(P)} \ \emph{For all $i$ and $j$ with $1<i<j\leq k$,
\begin{itemize}
\item[{\rm (a)}] $n(X_{i-1})\cap n(X_j)\subseteq n(X_{i-1})\cap
  n(X_i)$, and
\item[{\rm (b)}] $n(X_i)-n(X_{i-1})\subseteq n(X_j)-n(X_{i-1})$
  whenever $n(X_j)\cap n(X_{i-1})$ is nonempty.
\end{itemize}}
\end{quote}

Thus, to determine whether a transversal matroid $M$ with a given
presentation is a lattice path matroid, carry out the following steps.
\begin{itemize}
\item[{\rm (1)}] Detect and delete the isthmuses.
\item[{\rm (2)}] Determine the connected components.
\item[{\rm (3)}] Find the maximal presentation for each connected
  component.
\item[{\rm (4)}] For each component, find the classes defined above
  relative to the maximal presentation.
\item[{\rm (5)}] For each component, determine whether there is a
  linear order of these classes that satisfies property (P).
\item[{\rm (6)}] If there is such a linear order of these classes for
  each component, then use the criterion in Theorem~\ref{thm:charint}
  to determine whether, with respect to any corresponding linear order
  of a component, the intervals in the maximal presentation of that
  component are those of a maximal presentation of a lattice path
  matroid.
\end{itemize}
If, in step (5), there is no suitable order for some connected
component, then $M$ is not a lattice path matroid. If there is such an
order for each connected component, then $M$ is a lattice path matroid
if and only if step (6) yields only positive results.  Each of these
steps can be done in polynomial time in the size of the ground set, so
we get the following result.

\begin{thm}
  Whether a transversal matroid is a lattice path matroid can be
  determined from any presentation in polynomial time in the size of
  the ground set.
\end{thm}

\subsection{A Geometric Description of Lattice Path
  Matroids}\label{ssec:geom}

Brylawski~\cite{affine} (see also~\cite[Proposition 12.2.26]{ox}) gave
a geometric description of arbitrary transversal matroids.  This
section applies his result to lattice path matroids.

Let $M$ be a transversal matroid on the set $\{x_1,x_2,\ldots,x_k\}$
with presentation $(A_1,A_2,\ldots,A_r)$.  Brylawski showed that $M$
can be realized geometrically as follows.  Start with the free matroid
$M_0$ on a set $\{e_1,e_2,\ldots,e_r\}$ disjoint from $E(M)$.  For $i$
from $1$ to $k$, form $M_i$ from $M_{i-1}$ by taking the principal
extension of $M_{i-1}$ defined by the flat
$\cl_{M_{i-1}}(\{e_j\,:\,x_i\in A_j\})$, with the element added being
$x_i$.  The matroid $M$ is $M_k\del \{e_1,e_2,\ldots,e_r\}$.
Thus, a rank-$r$ matroid is transversal if and only if it can be
realized by placing the elements freely on the faces of the
$r$-simplex.

The next theorem, which is illustrated in Figure~\ref{draw}, shows how
lattice path matroids can be constructed by successively adding
isthmuses and loops, and by taking principal extensions by certain
flats.  To motivate this result, consider a lattice path matroid
$M[P,Q]$ that has rank $r$ and nullity $m$ in which $m+r$ is neither a
loop nor an isthmus. Let $l$ be the length of the longest final
segment of North steps in $P$. By Theorem~\ref{thm:max}, the sets of
the maximal presentation of $M[P,Q]$ that contain $m+r$ are the last
$l$ (those arising from $N_{r-l+1},\ldots,N_r$). By Brylawski's
result, $m+r$ is added freely to the flat spanned by $e_{r-l+1},
\ldots, e_r$ in the notation above; note that this flat is also
spanned by the last $l$ elements of $[m+r-1]$, since they are
independent in $M[P,Q]\del(m+r)$. Thus, we have the following result.

\begin{thm}\label{thm:extend}
  A matroid $M$ is a lattice path matroid if and only if the ground
  set can be written as $\{x_1,x_2,\ldots,x_k\}$ so that each
  restriction $M_i:=M|\{x_1,x_2,\ldots,x_i\}$ is formed from $M_{i-1}$
  by either
\begin{itemize}
\item[(i)] adding $x_i$ as an isthmus,
\item[(ii)] adding $x_i$ as a loop, or
\item[(iii)] adding $x_i$ via the principal extension of $M_{i-1}$
  generated by the closure of an independent set of the form $\{x_{h},
  x_{h+1},\ldots,x_{i-1}\}$ for some $h$ with $h<i$.
\end{itemize}
\end{thm}

\subsection{Relation to Other Classes of Transversal Matroids}\label{ssec:ftm}

We have seen that the class of lattice path matroids is closed under
taking both minors and duals.  While~\cite{omer} develops a
dual-closed, minor-closed class of transversal matroids that properly
contains $\mathcal{L}$, and while there are infinitely many
dual-closed, minor-closed classes contained in $\mathcal{L}$ (see
Sections~\ref{sec:gcat} and~\ref{sec:notch} for two such classes), few
other known classes of transversal matroids are either dual-closed or
minor-closed.  In this section, we make some remarks about two
important classes of transversal matroids, each of which has one of
these properties.

Fundamental transversal matroids (called principal transversal
matroids in~\cite{affine}) were introduced by Bondy and
Welsh~\cite{bw} and they play an important role in the study of
transversal matroids.  A transversal matroid $M$ is a
\emph{fundamental transversal matroid} if it can be represented on the
simplex with an element of $M$ at each vertex of the simplex.  Thus,
transversal matroids are the restrictions of fundamental transversal
matroids.  While the class $\mathcal{F}$ of fundamental transversal
matroids is closed under neither deletion nor contraction, it is
well-known and not hard to prove that $\mathcal{F}$ is dual-closed.
The class $\mathcal{F}$ is much larger than $\mathcal{L}$:
Brylawski~\cite{affine} showed that there are on the order of
$c^{n^2}$ simple fundamental transversal matroids on $n$ elements, for
some constant $c$; in contrast, $4^n$ is an upper bound on the number
of lattice path matroids on $n$ elements since there are $4^n$ pairs
of paths of length $n$ (see~\cite{lpm1} for a formula for the number of
connected lattice path matroids). Both $\mathcal{F}$ and $\mathcal{L}$ contain
all transversal matroids of rank two.  However, a fundamental
transversal matroid of rank three or more cannot have a pair of
disjoint connected hyperplanes, but such hyperplanes can occur in
lattice path matroids, such as the matroid $P_n=T_n(U_{n-1,n}\oplus
U_{n-1,n})$ of Theorem~\ref{thm:oxprow}.  On the other hand, the
number of connected hyperplanes of a fundamental transversal matroid,
such as the $n$-whirl $\mathcal{W}^n$, can exceed two (see
Corollary~\ref{cor:finterval}).

Let us call a matroid \emph{bitransversal} if both the matroid and its
dual are transversal.  It is easy to prove that the class of
bitransversal matroids is closed under direct sums, free extensions,
and free coextensions. Hence by starting with the union of the classes
$\mathcal{L}$ and $\mathcal{F}$, and using these three operations, we
can construct a larger class of bitransversal matroids; let
$\mathcal{LF}$ denote this class. For instance, the free extension
$(P_n \oplus \mathcal{W}^n) +e$ of $P_n \oplus \mathcal{W}^n$ is in
$\mathcal{LF}$ but not in $\mathcal{L}\cup \mathcal{F}$.  There are
bitransversal matroids, such as the identically self-dual matroids
of~\cite[Section 4]{bw}, that are not in $\mathcal{LF}$.  The problem
of characterizing all bitransversal matroids, which was posed by
Welsh, currently remains open (see~\cite[Problem 14.7.4]{ox}).

Bicircular matroids~\cite{bic} form another important class of
transversal matroids.  The notion of a bicircular matroid we consider
is a mild extension of that in~\cite{bic} (as originally defined,
bicircular matroids have no loops).  A transversal matroid $M$ is
\emph{bicircular} if it has a presentation $\mathcal{A}$ so that each
element of $M$ is in at most two sets in $\mathcal{A}$ (counting
multiplicity).  Thus, bicircular matroids are the transversal matroids
that have a representation on the simplex in which all nonloops are on
vertices or lines of the simplex.  It follows that minors of
bicircular matroids are bicircular.  On the other hand, the class of
bicircular matroids is not dual-closed: the prism (the matroid
$C_{4,2}$ of Figure~\ref{b22}) is bicircular, but its dual (the
matroid $B_{2,2}$ in the same figure) is not transversal.  Among the
matroids that are both bicircular and lattice path matroids are all
transversal matroids of rank two as well as iterated parallel
connections of rank-$2$ uniform matroids, $M_1:=U_{2,n_1}$ and
$M_i:=P(M_{i-1},U_{2,n_i})$, where the basepoint used to construct
$M_i$ is not in $M_{i-2}$.  A bicircular matroid, unlike a lattice
path matroid, can have more than two connected hyperplanes.  Also,
while most uniform matroids are not bicircular (for instance,
$U_{3,n}$ is bicircular if and only if $n\leq 6$), all uniform
matroids are in $\mathcal{L}$.  Thus, the class of bicircular matroids
differs significantly from $\mathcal{L}$ in all ranks greater than
two.

\section{Higher Connectivity}\label{sec:high}

In this section, we show how to find the connectivity $\lambda(M)$ of
a lattice path matroid in a simple way from the path presentation of
$M$.  We also show that at least one exact $\lambda(M)$-separation of
$M$ is given by a fundamental flat and its complement. We start by
recalling the relevant definitions; for more information on higher
connectivity, see~\cite[Chapter 8]{ox}.

For a positive integer $k$, a \emph{$k$-separation} of a matroid $M$ is
a partition of the ground set into two sets $X$ and $Y$, each with at
least $k$ elements, such that the inequality $r(X)+r(Y)\leq r(M)+k-1$
holds. A $k$-separation for which the equality $r(X)+r(Y)=r(M)+k-1$
holds is an \emph{exact $k$-separation}. The \emph{connectivity}, or
\emph{Tutte connectivity}, $\lambda(M)$ of $M$ is the least positive
integer $k$ such that $M$ has a $k$-separation; if there is no such
$k$, then $\lambda(M)$ is taken to be $\infty$. The connectivity of
uniform matroids is well known (see~\cite[Corollary 8.1.8]{ox}), so we
consider only lattice path matroids that are not uniform.  Also, as
justified by Theorem~\ref{thm:conncrit}, we focus exclusively on
lattice path matroids that are connected.

\begin{figure}
\begin{center}
  \epsfxsize 4.25truein \epsffile{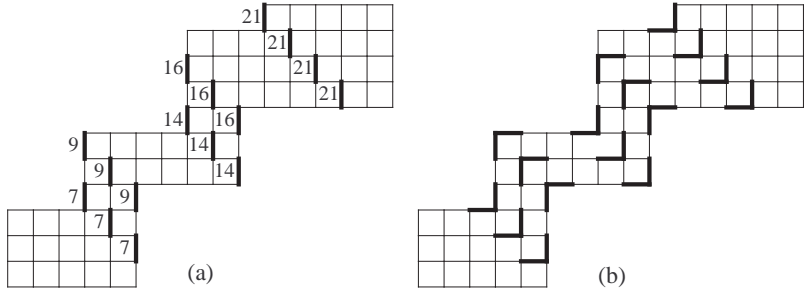}
\end{center}
\caption{(a) A lattice path matroid $M$ with $k_M=3$. (b) A
  pairing that shows the equality $k_M=k_{M^*}$. }\label{lambda}
\end{figure}

Let $M$ be a connected lattice path matroid, say $M[P,Q]$, that is not
uniform.  Let the integer $k_M$ be defined as follows:
$$
k_M:=\min\{|n(j)|\,:\, P \text{ has a } NE \text{ corner at } j
\text{ or } Q \text{ has an } EN \text{ corner at } j-1 \}.
$$
Figure~\ref{lambda} (a) illustrates a lattice path matroid $M$ in
which the relevant values of $j$ are $7,9,14,16$ (for which $|n(j)|$
is $3$) and $21$ (for which $|n(j)|$ is $4$), so $k_M$ is $3$.  The
main result of this section, Theorem~\ref{thm:conthm}, is that the
connectivity $\lambda(M)$ of $M$ is $k_M$.  Several lemmas enter into
the proof of this result. The first lemma reflects the equality
$\lambda(M)=\lambda(M^*)$ that holds for any matroid.

\begin{lemma}\label{lem:lambdadual}
  The number $k_M$ is invariant under duality, that is, $k_M=k_{M^*}$.
\end{lemma}

\begin{proof}
  Recall that the lattice path diagram for the dual of $M[P,Q]$ is
  obtained by reflecting the lattice path diagram for $M[P,Q]$ about
  the line $y=x$ (Figure~\ref{dual}).  Equivalently, the dual of
  $M[P,Q]$ is $M[Q',P']$ where $P'$ and $Q'$ are obtained from $P$ and
  $Q$ by switching East and North steps.  Let $n$ and $n'$ be the
  incidence functions of the standard presentations of $M[P,Q]$ and
  $M[Q',P']$, respectively.  Note that $P$ has a $NE$ corner at $j$ if
  and only if $P'$ has an $EN$ corner at $j$; also, $Q$ has an $EN$
  corner at $j-1$ if and only if $Q'$ has a $NE$ corner at at $j-1$.
  Thus, the lemma follows once we show the following statements: if
  $Q$ has an $EN$ corner at $j-1$, then $|n(j)|=|n'(j-1)|$; if $P$ has
  a $NE$ corner at $j$, then $|n(j)|=|n'(j+1)|$.  These assertions
  hold since we can pair off the relevant East and North steps that
  share a lattice point, as suggested in Figure~\ref{lambda} (b).
\end{proof}

Recall that in a matroid of connectivity at least $n$ with at least
$2(n-1)$ elements, circuits and cocircuits have at least $n$
elements~\cite[Proposition 8.1.6]{ox}.  The next lemma will be used to
show that circuits and cocircuits of a lattice path matroid $M$ have
at least $k_M$ elements.

\begin{lemma}\label{lem:nbd}
  Every element of $M$ is in at least $k_M-1$ sets in the maximal
  presentation of $M$.
\end{lemma}

\begin{proof}
  Let $M$ be $M[P,Q]$ and let $n$ and $n'$ be the incidence functions
  of its standard and maximal presentations. Let steps $q_N+1$,
  $q_E-1$, and $p_E+1$ be, respectively, the first East step of $Q$,
  the last North step of $Q$, and the first North step of $P$.  By the
  symmetry given by the order-reversing isomorphism of $M[P,Q]$ onto
  $M[Q^\rho,P^\rho]$ (see Theorem~\ref{thm:shape2}), it suffices to
  prove (a) if $i\leq q_N$, then $|n'(i)|\geq k_M-1$ (b) if $q_N < i <
  q_E$, then $|n(i)|\geq k_M-1$ and (c) if $q_E\leq i\leq p_E$, then
  $|n(i)|=r(M)$.  Theorem~\ref{thm:max} and the observation that $q_N$
  is at least $k_M-1$ prove part (a).  Part (c) is trivial.  The proof
  of part (b) uses the following easily-verified statements.
  \begin{itemize}
  \item[(i)] If the $j$-th and $(j+1)$-st steps of $Q$ are East, then
    $n(j+1)$ is either $n(j)$ or $n(j)-\min\bigl(n(j)\bigr)$, so we
    have $|n(j)|-1\leq |n(j+1)|\leq |n(j)|$.
  \item[(ii)] If the $j$-th step of $Q$ is North, then $n(j-1)$ is
    either $X$ or $X-\min(X)$ where $X$ is $\{h-1\,:\,h\in n(j)\}$, so
    $|n(j)|-1 \leq |n(j-1)| \leq |n(j)|$.
  \end{itemize}
  First assume that steps $i,i+1,\ldots,h$ of $Q$ are East and that
  step $h+1$ is North.  Thus, $Q$ has an $EN$ corner at $h$.
  Statements (i) and (ii) give the inequalities
  $$|n(i)|\geq |n(i+1)|\geq \cdots \geq |n(h)| \geq |n(h+1)|-1.$$
  Since $|n(h+1)|\geq k_M$, we have $|n(i)|\geq k_M-1$.  Finally, if
  the $i$-th step of $Q$ is North, a similar application of statement
  (ii) completes the proof of part (b).
\end{proof}

From Lemmas~\ref{lem:circform} and~\ref{lem:nbd}, the rank of any
circuit of $M$ is at least $k_M-1$.  The next lemma follows from this
observation and Lemma~\ref{lem:lambdadual}.  The generalized Catalan
matroid $M[(NE)^2]$ shows that $M$ can have circuits of rank $k_M-1$.

\begin{lemma}\label{lem:indlem}
  Any set of $k_M-1$ element of $[m+r]$ is independent in both $M$ and
  $M^*$. Circuits of $M$ have at least $k_M$ elements, as do
  circuits of $M^*$.
\end{lemma}

We now prove that $k_M$ is the connectivity of the lattice path
matroid $M$.

\begin{thm}\label{thm:conthm}
  Let $M$ be a connected lattice path matroid of rank $r$ and nullity
  $m$, say $M[P,Q]$, that is not uniform. The connectivity
  $\lambda(M)$ of $M$ is $k_M$, where $k_M$ is
  $$
  \min\{|n(j)|\,:\, P \text{ has a } NE \text{ corner at } j \text{
    or } Q \text{ has an } EN \text{ corner at } j-1 \}.
  $$
  Furthermore, at least one exact $k_M$-separation of $M$ consists
  of some fundamental flat and its complement.
\end{thm}

\begin{proof}
  We first show that $M$ has an exact $k_M$-separation that consists
  of a fundamental flat and its complement. Assume first that $k_M$ is
  $|n(j)|$ where $P$ has a $NE$ corner at $j$. Let $X$ and $Y$ be
  $[j]$ and $[j+1,m+r]$, respectively. Thus, $Y$ is a fundamental flat
  of $M$.  Note that both $X$ and $Y$ have at least $|n(j)|$ elements.
  It follows from the path presentations of restrictions given in
  Corollary~\ref{cor:restint} that $r(X)$ is $r(M)-r(Y)+|n(j)|-1$,
  that is, $r(X)+r(Y)=r(M)+k_M-1$, so $X,Y$ is an exact
  $k_M$-separation of $M$.  Similarly, if $k_M$ is $|n(j)|$ where $Q$
  has an $EN$ corner at $j-1$, then $[j-1]$ and $[j,m+r]$ give an
  exact $k_M$-separation of $M$.
  
  It remains to show that $M$ has no $h$-separation for any positive
  integer $h$ less than $k_M$. Let $h$ be such an integer and assume
  $X$ and $Y$ partition $[m+r]$, where both $X$ and $Y$ have at least
  $h$ elements.  We need to prove the inequality
  \begin{equation}\label{eqn:conngoal} 
   r(X)+r(Y)\geq r(M)+h. 
  \end{equation}
  If an element $y$ in $X$ is in the closure of $Y$, and if $X$ has
  more than $h$ elements, then we have $|X-y|\geq h$, $|Y\cup y|\geq
  h$, and $r(X)+r(Y)\geq r(X-y)+r(Y\cup y)$. Thus, it suffices to
  prove inequality~(\ref{eqn:conngoal}) when
   $|X|$ is $h$ or $Y$ is a nontrivial flat of $M$. By
  Lemma~\ref{lem:indlem}, each nontrivial connected component of the
  restriction $M|Y$ to a flat $Y$ of $M$ has more than $h$ elements;
  with an argument similar to the one above, it follows that if $Y$ is
  a nontrivial flat of $M$, then we may assume $Y$ is connected.
  
  Assume $|X|$ is $h$.  By Lemma~\ref{lem:indlem}, $X$ is an
  independent set that does not contain a cocircuit, so $Y$ spans $M$.
  Thus, $r(X)+r(Y)$ is $r(M)+h$.
  
  Now assume $Y$ is a nontrivial connected flat of $M$. If $Y$ is a
  fundamental flat, then inequality~(\ref{eqn:conngoal}) follows as in
  the first paragraph.  If $Y$ is not a fundamental flat, then, by
  Theorem~\ref{thm:cflats}, $Y$ is the intersection of two
  incomparable fundamental flats, say $Y\cup A$ and $Y\cup B$ where
  $A$ and $B$ partition $X$.  We may assume $1$ is in $A$, so $m+r$ is
  in $B$. Since $A\cup Y$ is a fundamental flat and $B$ is the
  complement of $A\cup Y$, we have $r(A\cup Y) + r(B) \geq
  r(M)+k_M-1$.  Thus, since $k_M$ exceeds $h$, to prove
  inequality~(\ref{eqn:conngoal}), it suffices to prove $r(X)+r(Y)\geq
  r(A\cup Y) + r(B)$, that is,
  \begin{equation}\label{eqn:conngoala}
    r(A\cup B)+r(Y)\geq r(A\cup Y) + r(B).
  \end{equation}
  Observe that $r(A\cup B)$ is $|n(A\cup B)|$; the inequality $r(A\cup
  B) \leq |n(A\cup B)|$ is obvious and the inequality $r(A\cup B) \geq
  |n(A\cup B)|$ follows by matching each set $N_i$, for $i$ in $n(A)$,
  with its first element, which must be in $A$, and each set in $N_j$,
  for $j$ in $n(B)-n(A)$, with its last element, which must be in $B$.
  A similar argument gives the equality $r(B) = |n(B)|$.  From
  Theorem~\ref{thm:finterval}, we also have $r(A\cup Y) = |n(A\cup
  Y)|$ and $r(Y) = |n(Y)|$. Thus, inequality~(\ref{eqn:conngoala}) is
  equivalent to
  \begin{equation}\label{eqn:conngoal1}
    |n(A\cup B)|+|n(Y)|\geq |n(A\cup Y)| + |n(B)|.
  \end{equation}
  Note that $|n(A\cup B)|$ is $|n(A)|+|n(B)|-|n(A)\cap n(B)|$.
  Substituting this and the analogous formula for $|n(A\cup Y)|$ into
  inequality~(\ref{eqn:conngoal1}) and simplifying gives that this
  inequality is equivalent to the inequality $|n(A)\cap n(Y)|\geq
  |n(A)\cap n(B)|$, which clearly holds. Thus,
  inequality~(\ref{eqn:conngoal}) holds, as needed to complete the
  proof.
\end{proof}

As the matroid $E_3$ of Figure~\ref{small} shows, not every exact
$k_M$-separation of a lattice path matroid $M$ has a fundamental flat
as one of the sets.

\section{Notch Matroids and their Excluded Minors}\label{sec:notch}

There are infinitely many minor-closed, dual-closed classes of
transversal matroids within the class of lattice path matroids.  One
way to define such classes is to impose certain requirements on the
bounding paths; for example, the lower bounding path of a generalized
Catalan matroid must have the form $E^mN^r$.  In this section we
introduce the minor-closed, dual-closed class of notch matroids, which
is defined by special forms for the bottom bounding path.  We relate
notch matroids to generalized Catalan matroids via circuit-hyperplane
relaxations.  The main result is the characterization of notch
matroids by excluded minors.  We include some remarks on the excluded
minors for lattice path matroids.

\begin{dfn}
  A \emph{notch matroid} is, up to isomorphism, a lattice path matroid
  of the form $M[E^mN^r,Q]$ or $M[E^{m-1}NEN^{r-1},Q]$.
\end{dfn}

\begin{figure}
\begin{center}
  \epsfxsize 4.0truein \epsffile{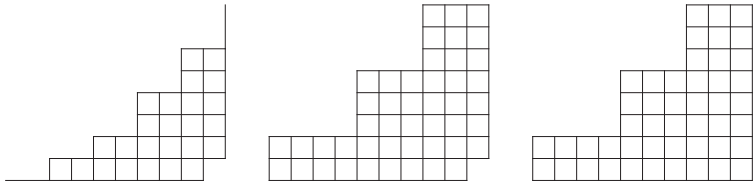}
\end{center}
\caption{Lattice path presentations of three notch
  matroids.}\label{notshp}
\end{figure}

\begin{figure}
\begin{center}
  \epsfxsize 4.0truein \epsffile{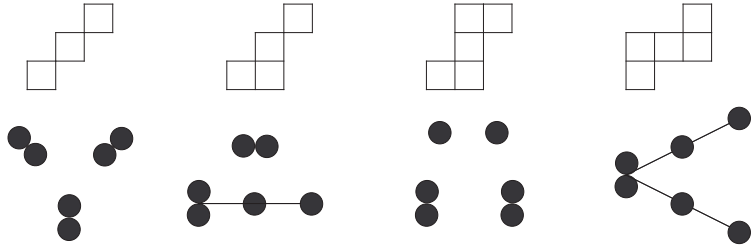}
\end{center}
\caption{Path presentations and geometric representations of
  $U_{1,2}\oplus U_{1,2}\oplus U_{1,2}$, $T_2(U_{1,2}\oplus
  U_{1,1}\oplus U_{1,1})\oplus U_{1,2}$, $D_3$, and
  $E_3$.}\label{small}
\end{figure}

As Figure~\ref{notshp} illustrates, notch matroids are either in
$\mathcal{C}$ or their lattice path presentations differ from those of
generalized Catalan matroids by the ``notch'' in the lower right
corner.  It follows from the lattice path descriptions of minors and
duals, along with Theorem~\ref{thm:shape2}, that the class
$\mathcal{N}$ of notch matroids is minor-closed and dual-closed.  Note
that $\mathcal{N}$, like its subclass $\mathcal{C}$, is not closed
under direct sums. In contrast to $\mathcal{C}$, the class
$\mathcal{N}$ is not closed under any of the following operations, as
can be seen from the matroid $D_3$ of Figure~\ref{small}: free
extension, truncation, and the dual operations.  The first lemma gives
a basic property that $\mathcal{N}$ shares with $\mathcal{C}$.

\begin{lemma}\label{lem:loopisth}
  Adding loops and isthmuses to a notch matroid yields a notch
  matroid.
\end{lemma}

\begin{figure}
\begin{center}
  \epsfxsize 3.0truein \epsffile{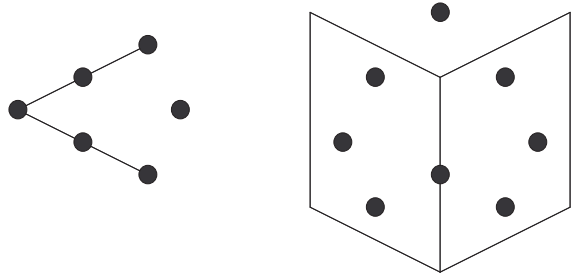}
\end{center}
\caption{The matroids $A_3$ and $A_4$.}\label{a3}
\end{figure}

Note that a connected notch matroid either is in $\mathcal{C}$ or has
a circuit-hyperplane relaxation in $\mathcal{C}$.  Not every matroid
that has a circuit-hyperplane relaxation in $\mathcal{C}$ is a notch
matroid; for instance, the matroids $A_3$ and $A_4$ of Figure~\ref{a3}
each have two circuit-hyperplane relaxations that are in
$\mathcal{C}$, yet neither is a lattice path matroid since condition
(ii) of Theorem~\ref{thm:lpmchar} fails.  However, we have the
following result.

\begin{thm}\label{thm:chrelax}
  A connected matroid in $\mathcal{L}-\mathcal{C}$ is a notch matroid
  if and only if it has a circuit-hyperplane.  Relaxing any
  circuit-hyperplane of a lattice path matroid yields a generalized
  Catalan matroid.
\end{thm}

\begin{proof}
  The last $r$ elements of a connected notch matroid
  $M[E^{m-1}NEN^{r-1},Q]$ obviously form a circuit-hyperplane.  For
  the converse, assume that $H$ is a circuit-hyperplane of the
  rank-$r$, nullity-$m$ matroid $M=M[P,Q]$.  Since $H$ is an
  $r$-circuit of $M$, by Theorem~\ref{thm:finterval} the set $n(H)$ is
  an interval of $r-1$ elements in $[r]$; we may assume that $n(H)$ is
  $[2,r]$. Since $H$ is a flat, $H$ is an interval of $r$ elements in
  the ground set $[m+r]$ of $M$, so $[m+r]$ consists of an initial
  interval, the interval $H$, and a final interval $Y$.  Since $H$ is
  a hyperplane, $Y$ must be empty, so $H$ consists of the last $r$
  elements of $[m+r]$. From these conclusions, it is immediate that
  $M$ is a notch matroid.  The last assertion follows from part (iii)
  of Corollary~\ref{cor:finterval}.
\end{proof}

Similar ideas yield the following result.

\begin{lemma}\label{lem:relax}
  Let $M'$ be $M[Q]$, a connected rank-$r$, nullity-$m$ matroid in
  $\mathcal{C}$.  If the basis $B$ of $M'$ is mapped onto the final
  segment $[m+1,m+r]$ by some automorphism of $M'$, then there is a
  unique matroid $M$ in which $B$ is a circuit-hyperplane and from
  which $M'$ is obtained by relaxing $B$.  Furthermore, $M$ is in
  $\mathcal{N}$.
\end{lemma}

The following two lemmas will be used heavily in the proof of the
excluded-minor characterization of $\mathcal{N}$.

\begin{lemma}\label{lem:incomplemma}
  If $X$ and $Y$ are nontrivial incomparable connected flats of a notch
  matroid $M$ that has no isthmuses, then either $X$ or $Y$ is a
  circuit-hyperplane.
\end{lemma}

\begin{proof}
  The incomparable flats $X$ and $Y$ show that $M$ is not in
  $\mathcal{C}$, so $M$ has a circuit-hyperplane, say $H$. Either $X$
  or $Y$ must be $H$ since $H$ cannot properly contain either $X$ or
  $Y$ and, by part (iii) of Corollary~\ref{cor:finterval}, nontrivial
  connected flats that are not contained in $H$ are comparable.
\end{proof} 

\begin{lemma}\label{lem:no3incomp}
  Three nontrivial connected flats $X$, $Y$, and $Z$ of a notch
  matroid $M$ cannot be mutually incomparable.
\end{lemma}

\begin{proof}
  We may assume that $M$ has no isthmuses and that $X$ and $Y$ are
  incomparable. From Lemma~\ref{lem:incomplemma}, either $X$ or $Y$,
  say $X$, is a circuit-hyperplane of $M$.  Part (iii) of
  Corollary~\ref{cor:finterval} implies that $Y$ and $Z$ are
  comparable.
\end{proof}

We turn to the excluded-minor characterization of $\mathcal{N}$.  Let
$ex(\mathcal{N})$ and $ex(\mathcal{L})$ denote the sets of excluded
minors for $\mathcal{N}$ and $\mathcal{L}$, respectively.  We first
discuss the matroids in $ex(\mathcal{N})$ that are not lattice path
matroids and so are in $ex(\mathcal{N})\cap ex(\mathcal{L})$.  In each
case, we show that the matroids are not in $\mathcal{L}$; it is easy to 
check that all their proper minors are in $\mathcal{N}$, so we omit
this part.

\begin{figure}
\begin{center}
  \epsfxsize 2.8truein \epsffile{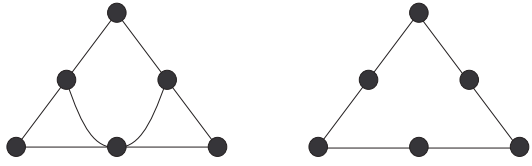}
\end{center}
\caption{The $3$-wheel $\mathcal{W}_3$ and the $3$-whirl
  $\mathcal{W}^3$.}\label{ww}
\end{figure}

Among the self-dual matroids in $ex(\mathcal{N})\cap ex(\mathcal{L})$
are the 3-wheel $\mathcal{W}_3$ and the 3-whirl $\mathcal{W}^3$, which
are shown in Figure~\ref{ww}. Since all $3$-point lines of
$\mathcal{W}_3$ and $\mathcal{W}^3$ are fundamental flats, condition
(i) of Theorem~\ref{thm:lpmchar} fails, so $\mathcal{W}_3$ and
$\mathcal{W}^3$ are not in $\mathcal{L}$.

For $n\geq 3$, let $A_n$ be the rank-$n$ paving matroid with only two
nontrivial hyperplanes, $\{x,a_2,a_3,\ldots,a_n\}$ and
$\{x,b_2,b_3,\ldots,b_n\}$, and with only one point, $y$, in neither
circuit-hyperplane (Figure~\ref{a3}).  The two circuit-hyperplanes
violate condition (ii) of Theorem~\ref{thm:lpmchar}, so $A_n$ is not
in $\mathcal{L}$.  Note that $A_n$ is self-dual.

\begin{figure}
\begin{center}
  \epsfxsize 4.0truein \epsffile{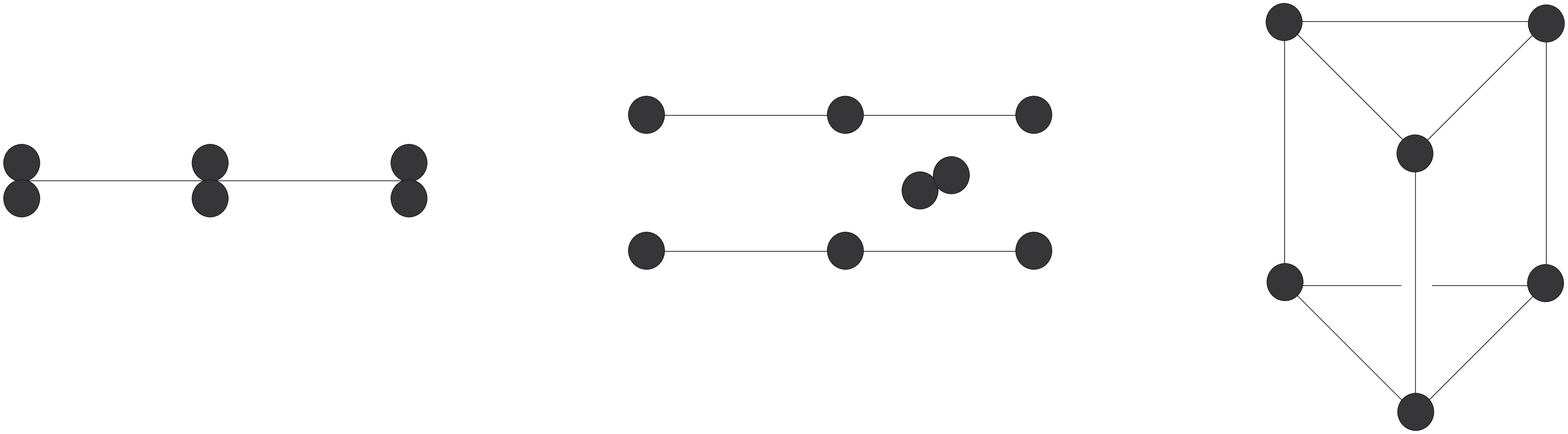}
\end{center}
\caption{The matroids $B_{2,2}$, $B_{3,2}$, and
  $C_{4,2}$.}\label{b22}
\end{figure}

We next consider two doubly-indexed families in $ex(\mathcal{N})\cap
ex(\mathcal{L})$ that are related by duality; three of these matroids
are shown in Figure~\ref{b22}.  Let $n$ and $k$ be integers with
$2\leq k\leq n$.  Let $B_{n,k}$ be the truncation $T_n(U_{n-1,n}\oplus
U_{n-1,n}\oplus U_{k-1,k})$ to rank $n$ of the direct sum of two
$n$-circuits and a $k$-circuit.  The three disjoint circuits are
fundamental flats of $B_{n,k}$, so condition (i) of
Theorem~\ref{thm:lpmchar} shows that $B_{n,k}$ is not in
$\mathcal{L}$.  The dual $C_{n+k,k}$ of $B_{n,k}$ is the rank-$(n+k)$
paving matroid $C_{n+k,k}$ for which the ground set can be partitioned
into sets $X,Y,Z$ with $|X| = |Y| = n$ and $|Z| = k$ so that the only
nontrivial hyperplanes are $X\cup Y$, $X\cup Z$, and $Y\cup Z$.

\begin{figure}
\begin{center}
  \epsfxsize 3.0truein \epsffile{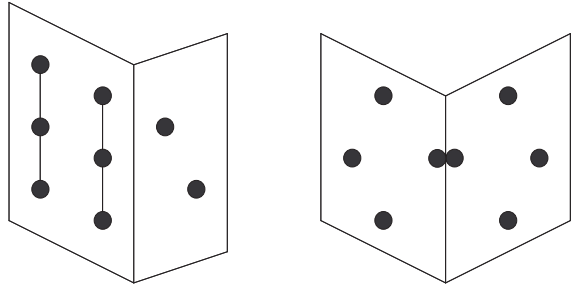}
\end{center}
\caption{The matroids $D_4$ and $E_4$.}\label{d4}
\end{figure}

The remaining matroids in $ex(\mathcal{N})\cap ex(\mathcal{L})$, two
of which are shown in Figure~\ref{d4}, form two infinite families that
are related by duality. Recall that $M+y$ denotes the free extension
of $M$ by the point $y$. For $n\geq 3$, let $D_n$ be the rank-$n$
matroid
$$\bigl(T_{n-1}(U_{n-2,n-1}\oplus U_{n-2,n-1})\oplus
U_{1,1}\bigr)+y.$$
That $D_n$ is not in $\mathcal{L}$ for $n\geq 4$
follows since the two $(n-1)$-circuits, as well as their union, are
fundamental flats of $D_n$, contrary to condition (i) of
Theorem~\ref{thm:lpmchar}.  In the dual $E_n$ of $D_n$, the element
$y$ is parallel to an element $x$, and the deletion $E_n\del y$ is a
rank-$n$ paving matroid whose only nontrivial hyperplanes are two
circuit-hyperplanes that intersect in $x$.  (The matroids $D_3$ and
$E_3$, which are shown in Figure~\ref{small}, are lattice path
matroids.)

We have proven the easy part of the following theorem; the more
substantial part of this result follows from the excluded-minor
characterization of notch matroids, which is given in
Theorem~\ref{thm:notchchar}.  

\begin{thm}\label{thm:exlexn}
  The matroids in $ex(\mathcal{L})\cap ex(\mathcal{N})$ are:
\begin{itemize}
\item[($1$)] the three-wheel $\mathcal{W}_3$ and the three-whirl
  $\mathcal{W}^3$,
\item[($2$)] $A_n$ for $n\geq 3$,
\item[($3$)] $B_{n,k}$ and $C_{n+k,k}$ for $n$ and $k$ with $2\leq
  k\leq n$, and
\item[($4$)] $D_n$ and $E_n$ for $n\geq 4$.
\end{itemize}
\end{thm}

\begin{figure}
\begin{center}
  \epsfxsize 3.75truein \epsffile{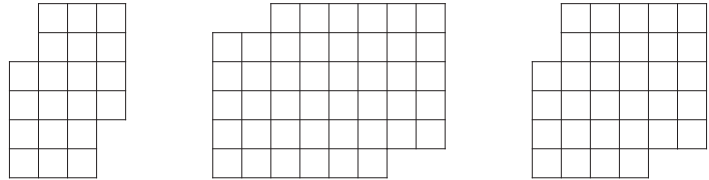}
\end{center}
\caption{Lattice path presentations of $F_6$, $G_6$, and
  $H_6$.}\label{large}
\end{figure}

We now turn to the excluded-minor characterization of notch matroids.
The excluded minors are those in Theorem~\ref{thm:exlexn} together
with the three types of lattice path matroids illustrated in
Figure~\ref{large} and the four  matroids in Figure~\ref{small}.

\begin{thm}\label{thm:notchchar}
  The excluded minors for the class of notch matroids are:
\begin{itemize}
\item[($1$)] $U_{1,2}\oplus U_{1,2}\oplus U_{1,2}$ and
  $T_2(U_{1,2}\oplus U_{1,1}\oplus U_{1,1})\oplus U_{1,2}$,
\item[($2$)] the three-wheel, $\mathcal{W}_3$, and the three-whirl,
  $\mathcal{W}^3$,
\item[($3$)] $A_n$ for $n\geq 3$,
\item[($4$)] $B_{n,k}$ and $C_{n+k,k}$ for $n$ and $k$ with $2\leq
  k\leq n$,
\item[($5$)] $D_n$ for $n\geq 3$,
\item[($6$)] $E_n$ for $n\geq 3$,
\item[($7$)] for $n\geq 4$, the rank-$n$ matroid $F_n:=
  T_n(U_{n-2,n-1}\oplus U_{n-2,n-1})$,
\item[($8$)] for $n\geq 2$, the rank-$n$ matroid $G_n:=
  T_n(U_{n-1,n+1}\oplus U_{n-1,n+1})$, and
\item[($9$)] for $n\geq 3$, the rank-$n$ matroid $H_n:=
  T_n(U_{n-2,n-1}\oplus U_{n-1,n+1})$.
\end{itemize}
\end{thm}

To make the proof of Theorem~\ref{thm:notchchar} less verbose, we will
use abbreviations such as the following: from Theorem~\ref{thm:notlpm}
applied to $M$, $X_1$, $X_2$, and $y$, we get $M\not\in\mathcal{L}$.
By this we mean that the matroid $M$ and the flats $X_1$ and $X_2$
satisfy the hypotheses of Theorem~\ref{thm:notlpm}, with the point $y$
showing the validity of the third condition.

\begin{proof}[Proof of Theorem~\ref{thm:notchchar}.]
  The remarks before Theorem~\ref{thm:exlexn} show that of the
  matroids in the theorem, only $D_3$, $E_3$, and those in items ($1$)
  and ($7$)--($9$) are in $\mathcal{L}$.  The presentations of these
  matroids, illustrated in Figures~\ref{large} and~\ref{small}, make
  it clear that they are not in $\mathcal{N}$.  It is easy to check
  that all proper minors of these matroids are in $\mathcal{N}$.  Note
  that $H_n$ is self-dual, and that $F_n$ and $G_{n-2}$ are dual to
  each other.
  
  The proof that Theorem~\ref{thm:notchchar} gives all excluded minors
  is intricate, so we first outline the argument. Part
  (\ref{thm:notchchar}.1) proves that the disconnected excluded minors
  are $U_{1,2}\oplus U_{1,2}\oplus U_{1,2}$, $T_2(U_{1,2}\oplus
  U_{1,1}\oplus U_{1,1})\oplus U_{1,2}$, $F_4$, $G_2$, and $H_3$. The
  rest of the proof revolves around three properties a connected
  excluded minor $M$ may have:
  \begin{itemize}
  \item[(a)] $r(X_1\cup X_2)<r(M)$ for some nontrivial incomparable
    connected flats $X_1,X_2$,
   \item[(b)] $M$ contains three mutually incomparable connected
     flats,
  \item[(c)] $M$ has no circuit-hyperplane.
  \end{itemize}
  In (\ref{thm:notchchar}.2), we show that if $M$ has property (a),
  then $M$ is $D_n$ for some $n\geq 3$.  Part (\ref{thm:notchchar}.3)
  gives a key property of all connected excluded minors. In
  (\ref{thm:notchchar}.4), we show that if property (b) but not (a)
  holds, then $M$ is one of the matroids in items ($2$) and ($4$).
  Part (\ref{thm:notchchar}.5) shows that if only property (c) holds,
  then $M$ is one of the matroids in items ($6$)--($9$). If none of
  the properties holds, then for any mutually incomparable connected
  flats $X_1,X_2,\ldots,X_k$, we have $k\leq 2$, and if $k$ is $2$,
  then at least one of $X_1$ or $X_2$ is a circuit-hyperplane. Since
  restrictions to proper subsets of circuit-hyperplanes are free, it
  follows that relaxing a circuit-hyperplane of such an excluded minor
  yields a matroid $M'$ in which the connected flats are linearly
  ordered by inclusion, that is, $M'$ is in $\mathcal{C}$.  The proof
  of Theorem~\ref{thm:notchchar} is completed in
  (\ref{thm:notchchar}.6) by showing that the only rank-$n$ excluded
  minor that has a circuit-hyperplane relaxation in $\mathcal{C}$ is
  $A_n$.
  
  Throughout the proof, $M$ denotes a rank-$n$ excluded minor for the
  class of notch matroids.  By Lemma~\ref{lem:loopisth}, $M$ has
  neither loops nor isthmuses.

\begin{quote}
  {\bf (\ref{thm:notchchar}.1) } If $M$ is disconnected, then $M$ is
  one of $U_{1,2}\oplus U_{1,2}\oplus U_{1,2}$, $T_2(U_{1,2}\oplus
  U_{1,1}\oplus U_{1,1})\oplus U_{1,2}$, $F_4$, $G_2$, and $H_3$.
\end{quote}

\begin{proof}[Proof of (\ref{thm:notchchar}.1).]
  Assume $M$ has at least three components. Each component has a
  circuit of two or more elements, so $M$ has $U_{1,2}\oplus
  U_{1,2}\oplus U_{1,2}$ as a minor, which is itself an excluded
  minor.  Thus, $M$ is $U_{1,2}\oplus U_{1,2}\oplus U_{1,2}$.
  
  Now assume $M$ has exactly two components, $M_1$ and $M_2$. Being
  proper minors of $M$, both $M_1$ and $M_2$ are notch matroids.
  Observe that if $r(M_i)\geq 2$, then, by Theorem~\ref{thm:spanning}
  and Corollary~\ref{cor:cinterval}, there is an element $x$ for which
  $M_i/x$ is connected.  Dually, if $\eta(M_i)\geq 2$, then $M_i\del
  y$ is connected for some $y$.
  
  Assume $M_1$ is $U_{1,2}$.  From lattice path presentations and from
  the statements $M_2\in\mathcal{N}$ and $U_{1,2}\oplus
  M_2\not\in\mathcal{N}$, it follows that $r(M_2)$ and $\eta(M_2)$ are
  both at least $2$. Similarly, if $M'_2$ is a connected minor of
  $M_2$ for which $r(M'_2)$ and $\eta(M'_2)$ are both $2$, then
  $U_{1,2}\oplus M'_2\not\in\mathcal{N}$. These observations, together
  with those in the last paragraph, imply that $r(M_2)$ and
  $\eta(M_2)$ are both $2$.  From lattice path presentations, we see
  that only two connected lattice path matroids have rank and nullity
  $2$, namely $U_{2,4}$ and $T_2(U_{1,2}\oplus U_{1,1}\oplus
  U_{1,1})$, so $M$ is either $H_3$ or $T_2(U_{1,2}\oplus
  U_{1,1}\oplus U_{1,1})\oplus U_{1,2}$.
  
  Now assume $M_1=U_{1,k}$ with $k\geq 3$.  Since
  $M\not\in\mathcal{N}$, the nullity of $M_2$ is at least $2$.
  Arguments like those in the last paragraph imply that $k$ is $3$,
  that $\eta(M_2)$ is $2$, and that $r(M_2)$ is $1$; therefore $M_2$
  is $U_{1,3}$, so $M$ is $G_2$.
  
  Finally, if $M_1$ and $M_2$ have rank $2$ or greater, then, by the
  same types of arguments, both $M_1$ and $M_2$ have rank $2$ and
  nullity $1$, so $M$ is $F_4$.
\end{proof}

From now on, we assume $M$ is connected.

\begin{quote}
  {\bf (\ref{thm:notchchar}.2) } If $M$ has nontrivial incomparable
  connected flats $X_1$ and $X_2$ with $r(X_1\cup X_2)<n$, then $M$ is
  $D_n$.
\end{quote}

\begin{proof}[Proof of (\ref{thm:notchchar}.2).]
  Choose such a pair of flats $X_1,X_2$ so that $r(X_1) + r(X_2)$ is
  as small as possible.  Lemma~\ref{lem:incomplemma} applied to
  $M|(X_1\cup X_2)$, $X_1$, and $X_2$ implies that either $X_1$ or
  $X_2$ is a circuit-hyperplane of $M|(X_1\cup X_2)$.

  Assume $M|(X_1\cup X_2)$ is disconnected. This disconnected notch
  matroid has neither loops nor isthmuses, so one component, say
  $X_1$, has rank $1$ and the other, $X_2$, has nullity $1$; thus,
  $X_1$ is a parallel class and $X_2$ is a circuit.  If $|X_1|>2$ and
  $y\in X_1$, then $M\del y$, $X_1-y$, and $X_2$ contradict
  Lemma~\ref{lem:incomplemma}.  If $|X_2|>2$ and $z\in X_2$, then
  $M/z$, $\cl_{M/z}(X_1)$, and $X_2-z$ contradict
  Lemma~\ref{lem:incomplemma}.  Thus, $|X_1| = |X_2| = 2$.  Since $M$
  has neither $B_{2,2}$ nor $U_{1,2}\oplus U_{1,2}\oplus U_{1,2}$ as a
  proper minor, $X_1$ and $X_2$ are the only nontrivial parallel
  classes of $M$.  Let $x$ and $y$ be in $E(M) - \cl(X_1\cup X_2)$.
  By Lemma~\ref{lem:incomplemma}, the rank-$1$ flats $\cl_{M/x}(X_1)$
  and $\cl_{M/x}(X_2)$ are hyperplanes of $M/x$, so $r(M)$ is $3$.  It
  follows that $M|(X_1\cup X_2\cup\{x,y\})$, and so $M$, is one of the
  excluded minors $T_2(U_{1,2}\oplus U_{1,1}\oplus U_{1,1})\oplus
  U_{1,2}$ or $D_3$; since $M$ is connected, $M$ is $D_3$.
  
  Now assume $M|(X_1\cup X_2)$ is connected.  We show that $M$ is
  $D_n$ by proving the following statements:
  \begin{itemize}
   \item[(i)] $M$ is simple,
   \item[(ii)] $X_1$ and $X_2$ are disjoint circuits, and $X_1\cup
     X_2$ is a flat of $M$,
   \item[(iii)] $E(M) - (X_1\cup X_2)$ contains only two elements, say $x$
     and $y$,
   \item[(iv)] the only nonspanning circuits of $M\del x,y$ are $X_1$ and
     $X_2$,
   \item[(v)] $|X_1| = |X_2|$, so both $X_1$ and $X_2$ are
     circuit-hyperplanes of $M\del x,y$, and
   \item[(vi)] the only circuits of $M$ that contain $x$ and $y$ are
     spanning circuits.
  \end{itemize}
  
  To prove statement (i), note that since $M|(X_1\cup X_2)$ is
  connected, and since $X_1$ and $X_2$ are incomparable flats, neither
  $X_1$ nor $X_2$ is a parallel class.  If elements $x$ and $y$ of $M$
  were parallel, then $M\del y$, $X_1-y$, and $X_2-y$ (which may be
  $X_1$ and $X_2$) would contradict Lemma~\ref{lem:incomplemma}.
  
  For statement (ii), we first show that both $M|X_1/x$ and $M|X_2/x$
  are connected for any $x$ in $X_1\cap X_2$.  If, say, $M|X_1/x$ were
  disconnected, then by Lemma~\ref{lem:smallerflats}, there would be
  nontrivial incomparable connected flats $A$ and $B$ of $M|X_1$ with
  $r(A)+r(B) = r(X_1)+1$.  Since $M$ is simple, $r(X_2)$ exceeds $1$,
  so the flats $A$ and $B$ of $M$ would contradict the choice of $X_1$
  and $X_2$ as minimizing the sum $r(X_1)+r(X_2)$.  Since $M|X_1/x$
  and $M|X_2/x$ are connected, $M/x$, $X_1-x$, and $X_2-x$ contradict
  Lemma~\ref{lem:incomplemma}.  Thus, $X_1$ and $X_2$ are disjoint.
  The connected notch matroids $M|X_1$ and $M|X_2$ have spanning
  circuits; this observation and the minimality of $M$ show that $X_1$
  and $X_2$ are circuits.  For any $x$ in $\cl(X_1\cup X_2) - (X_1\cup
  X_2)$, the deletion $M\del x$ is connected, so $M\del x$, $X_1$, and
  $X_2$ would violate Lemma~\ref{lem:incomplemma}. Thus, $\cl(X_1\cup
  X_2)$ is $X_1\cup X_2$, so statement (ii) holds.
  
  Let $y$ be in $E(M) - (X_1\cup X_2)$. The contraction $M/y$ has
  neither loops nor isthmuses.  By Lemma~\ref{lem:incomplemma}, at
  least one of $\cl_{M/y}(X_1)$ and $\cl_{M/y}(X_2)$ is a
  circuit-hyperplane of the notch matroid $M/y$, so $r(X_1\cup X_2)$
  is $n - 1$.  For $M\del y$, $X_1$, and $X_2$ to not contradict
  Lemma~\ref{lem:incomplemma}, $M\del y$ must have an isthmus.  From
  these conclusions, statement (iii) follows.
  
  Assume $C$ is a nonspanning circuit of $M\del x,y$ other than $X_1$
  and $X_2$.  Recall that either $X_1$ or $X_2$, say $X_1$, is a
  circuit-hyperplane of $M\del x,y$.  Thus, $X_1$ and $\cl(C)$ are
  incomparable and $X_1\cup C$ spans the flat $X_1\cup X_2$.  Let $z$
  be in the difference $X_2-C$ of circuits.  Note that $M\del z$ is
  connected.  That $M\del z$, $X_1$, and $\cl(C)-z$ contradict
  Lemma~\ref{lem:incomplemma} proves statement (iv).  Statement (v)
  follows since if $|X_2|<|X_1|$ and $z$ is in $X_1$, then $M/z$,
  $X_1-z$, and $X_2$ would contradict Lemma~\ref{lem:incomplemma}.
  
  From statements (i) and (v) we have $n\geq 4$.  Assume $x$ and $y$
  are in a nonspanning circuit $C$. At least one of $X_1$ and $X_2$ is
  not contained in $\cl(C)$, so we may assume that $X_1$ and $\cl(C)$
  are incomparable.  Let $z$ be in the difference $X_2-C$ of circuits.
  Note that $X_1\cup (X_2-z)$ is a connected hyperplane of $M\del z$
  since $n\geq 4$, so $M\del z$ is connected.
  Lemma~\ref{lem:incomplemma} applied to $M\del z$, $X_1$ and
  $\cl_{M\del z}(C)$ implies that $\cl_{M\del z}(C)$ must be a
  circuit-hyperplane of $M\del z$, so $\cl(C)$ is a hyperplane of $M$.
  Note that $\cl(C)$ is either $\cl_{M\del z}(C)$ or $\cl_{M\del
    z}(C)\cup z$, that is, either $C$ or $C\cup z$, so $|\cl(C)|\leq
  n+1$.  Thus, if $X_2\subseteq \cl(C)$, then $\cl(C)$ is
  $X_2\cup\{x,y\}$.  However, if $\cl(C)$ is $X_2\cup\{x,y\}$ and $w$
  is in $X_1$, then $M\del w$, $(X_1-w)\cup X_2$, $X_2\cup\{x,y\}$
  contradict Lemma~\ref{lem:incomplemma}. Therefore $X_2$ and $\cl(C)$
  are incomparable.  By switching $X_1$ and $X_2$ if necessary, we may
  assume $C\cap X_1\ne\emptyset$. Since $r(C)=n-1$, we have $r(C\cup
  X_1) = n$; however, there are at least two elements, say $a$ and
  $b$, in $X_2 - \bigl(\cl(C)\cup X_1\bigr)$, that is, in $X_2 -
  \cl(C)$, so by Theorem~\ref{thm:notlpm}, $M\del a$ is not a lattice
  path matroid, contrary to the minimality of $M$. Thus, statement
  (vi) holds, so $M$ is $D_n$.
\end{proof}

\begin{quote} {\bf (\ref{thm:notchchar}.3) } 
  If $X$ is a proper nontrivial connected flat of $M$ and the
  element $x$ of $X$ is not parallel to any element, then $X-x$ is a
  connected flat of $M/x$.
\end{quote}

\begin{proof}[Proof of (\ref{thm:notchchar}.3).]
  If $X-x$ were a disconnected flat of $M/x$, then, by
  Lemma~\ref{lem:smallerflats} applied to $M|X$, we would have
  $r(X_1\cup X_2)\leq r(X)<r(M)$ for some nontrivial incomparable
  connected flats $X_1, X_2$ of $M|X$.  Since $X_1$ and $X_2$ would
  also be flats of $M$, by (\ref{thm:notchchar}.2), $M$ would be
  $D_n$.  That $D_n$ has no such flat $X$ and element $x$ provides the
  contradiction that proves the result.
\end{proof}

\begin{quote} {\bf (\ref{thm:notchchar}.4) } If $M$ 
  has three mutually incomparable connected flats $X_1,X_2,X_3$, then
  $M$ is $\mathcal{W}_3$, $\mathcal{W}^3$, $B_{n,k}$, or $C_{n,k}$.
\end{quote}

\begin{proof}[Proof of (\ref{thm:notchchar}.4).]
  The minimality of $M$ and Lemma~\ref{lem:no3incomp} imply that the
  ground set of $M$ is $X_1\cup X_2\cup X_3$ and that any pair $x,y$
  of parallel elements can be in only one of $X_1,X_2,X_3$.  If an
  element $x$ were in $X_1\cap X_2\cap X_3$, then by
  (\ref{thm:notchchar}.3), $M/x$, $X_1-x$, $X_2-x$, and $X_3-x$ would
  contradict Lemma~\ref{lem:no3incomp}, so $X_1\cap X_2\cap X_3 =
  \emptyset$. Note that $M$ is not $D_n$, so we have $r(X_i\cup X_j) =
  n$ for $\{i,j\} \subset \{1,2,3\}$.
  
  First assume $X_1\cap X_2 = \emptyset$.  There are at least two
  points $x$ and $y$ in $X_2 - X_3$, so if $X_1\cap X_3$ were
  nonempty, then $M\del y$, $X_1$, and $X_3$ would contradict
  Theorem~\ref{thm:notlpm}. Thus, $X_1\cap X_3 = \emptyset$.
  Similarly $X_2\cap X_3 = \emptyset$.  The minimality of $M$ implies
  that $X_1$, $X_2$, and $X_3$ are circuits. Let $\{i,j,k\}$ be
  $\{1,2,3\}$. Since $r(X_i\cup X_j)$ is $n$, for any $x$ in $X_k$ the
  notch matroid $M\del x$ has no isthmuses; thus, from
  Lemma~\ref{lem:incomplemma}, either $X_i$ or $X_j$ is a
  circuit-hyperplane of $M\del x$ and so of $M$. It follows that at
  least two of $X_1,X_2,X_3$, say $X_1$ and $X_2$, are
  circuit-hyperplanes of $M$.  Let $|X_3|$ be $k$. Note that $M$ is
  $B_{n,k}$ if $X_1$, $X_2$, and $X_3$ are the only nonspanning
  circuits of $M$.  If $C$ were another nonspanning circuit, then for
  any $z$ in the difference $X_3-C$ of circuits, the flat
  $\cl_{M\del z}(C)$ would be contained in neither of the hyperplanes
  $X_1$ and $X_2$ of $M\del z$, contrary to part (iv) of
  Corollary~\ref{cor:finterval}.  Thus, $M$ is $B_{n,k}$.
  
  Now assume $X_i\cap X_j\ne\emptyset$ for all sets
  $\{i,j\}\subset\{1,2,3\}$.  We claim that $X_1$, $X_2$, and $X_3$
  are hyperplanes and the union of any two contains all but at most
  one point of $M$. Let $\{i,j,k\}$ be $\{1,2,3\}$ and let $x$ be in
  $X_i\cap X_j$. The equality $r(X_i\cup X_k)=n$ and
  Theorem~\ref{thm:notlpm} give the inequality $|E(M)-(X_i\cup
  X_k)|\leq 1$, so the second claim holds.  To see that $X_k$ is a
  hyperplane, note that Lemma~\ref{lem:no3incomp} applied to $M/x$,
  $\cl_{M/x}(X_k)$, $X_i -x$, and $X_j - x$ implies that there is a
  containment among at least two of these sets.  Of the two possible
  containments, we may assume $X_i-x \subseteq \cl_{M/x}(X_k)$.  Thus,
  $X_i\subseteq \cl(X_k\cup x)$.  This containment, the inequality
  $|E(M)-(X_i\cup X_k)|\leq 1$, and that $X_j$ is connected imply that
  $\cl(X_k\cup x)$ is $X_1\cup X_2\cup X_3$, so $X_k$ is a hyperplane
  of $M$.
 
  If $x$ and $y$ are in $X_i\cap X_j$, then $x$ is in the nontrivial
  connected hyperplanes $X_i-y$ and $X_j-y$ of the notch matroid
  $M/y$, so, by Theorem~\ref{thm:notlpm}, $E(M/y)$ is $(X_i-y)\cup
  (X_j-y)$. Thus, if $|X_i\cap X_j|\geq 2$, then $E(M) = X_i\cup X_j$.
  
  Assume $|X_1\cap X_2|$ is $1$.  Since $X_1$ is connected and at most
  one point of $X_1$ is in neither $X_1\cap X_2$ (one point) nor
  $X_1\cap X_3$ (a flat), there is one point in $X_1 - (X_2\cup X_3)$.
  Similarly, there is one point in $X_2 - (X_1\cup X_3)$.  These
  conclusions, and that in the last paragraph, give the equality
  $|X_1\cap X_3|=|X_2\cap X_3|=1$.  Therefore $X_1$, $X_2$, and $X_3$
  are $3$-point lines. It follows easily that $M$ is either
  $\mathcal{W}_3$ or $\mathcal{W}^3$.
  
  Assume $|X_i\cap X_j|\geq 2$ for $\{i,j,k\} = \{1,2,3\}$. Thus, $X_i
  = (X_i\cap X_j)\cup (X_i\cap X_k)$.  Let $x$ and $y$ be in $X_i\cap
  X_j$.  Lemma~\ref{lem:incomplemma} applied to $M/y$, $X_i - y$, and
  $X_j - y$ implies that either $X_i - y$ or $X_j-y$ is a
  circuit-hyperplane of $M/y$.  Since, in addition, $X_i$ and $X_j$
  are connected hyperplanes of $M$, either $X_i$ or $X_j$ is a
  circuit-hyperplane of $M$.  It follows that at least two
  hyperplanes, say $X_1$ and $X_2$, are circuit-hyperplanes of $M$.
  Assume $|X_1\cap X_2| = k$.  That $X_1$ and $X_2$ are
  circuit-hyperplanes of $M$ gives the equality $|X_1\cap X_3| = n-k =
  |X_2\cap X_3|$. To prove that $M$ is $C_{n,k}$, we need only show
  that the only proper nontrivial connected flat $X$ other than $X_1$
  and $X_2$ is $X_3$.  Clearly $X$ is incomparable to the
  circuit-hyperplanes $X_1$ and $X_2$. As we deduced for
  $X_1,X_2,X_3$, we get $X\cap X_1\cap X_2 = \emptyset$, so
  $X\subseteq X_3$. Since $X_1\cap X_3$ and $X_2\cap X_3$ are
  independent, both $X\cap X_1$ and $X\cap X_2$ are nonempty.  With
  this, the claim in the third paragraph shows that $X$ is a
  hyperplane.  Since $X\subseteq X_3$, it follows that $X$ is $X_3$,
  as needed.
\end{proof}

\begin{quote}
  {\bf (\ref{thm:notchchar}.5) } If $M$ has no
  circuit-hyperplane and is not $D_n$, then $M$ is one of $E_n$,
  $F_n$, $G_n$, or $H_n$.
\end{quote}

\begin{proof}[Proof of (\ref{thm:notchchar}.5).]
  Since $M$ is not a generalized Catalan matroid, there is a pair
  $X_1$, $X_2$ of incomparable connected flats.  Since $M$ is not
  $D_n$, part (\ref{thm:notchchar}.2) gives the equality $r(X_1\cup
  X_2) = n$ for any such pair of flats.
  
  Assume there were an element $x$ in $E(M) - (X_1\cup X_2)$.  Since
  $r(X_1\cup X_2)$ is $n$, the deletion $M\del x$ would have no
  isthmuses. Therefore either $X_1$ or $X_2$ would be a
  circuit-hyperplane of $M\del x$ and so of $M$.  Since $M$ has no
  circuit-hyperplane, the equality $E(M) = X_1\cup X_2$ follows.
  
  First assume $M$ has two incomparable connected flats $X_1$ and
  $X_2$ that are not disjoint.  We show that $M$ is $E_n$ by proving
  the following statements:
   \begin{itemize} 
   \item[(i)] each element in $X_1\cap X_2$ is parallel to another
     element of $M$,
   \item[(ii)] $X_1\cap X_2$ contains just two elements, say $x$ and
     $y$, and at least one of $X_1-x$ and $X_2-x$, say $X_1-x$, is a
     circuit-hyperplane of $M\del x$,
   \item[(iii)] $X_2 - x$ is a circuit,
   \item[(iv)] $|X_1|=|X_2|$, and
   \item[(v)] the nonspanning circuits of $M$ are $X_1-x$, $X_1-y$,
     $X_2-x$, $X_2-y$, and $\{x,y\}$.
   \end{itemize}
   Assume statement (i) failed for some $x$ in $X_1\cap X_2$.  From
   (\ref{thm:notchchar}.3) and Lemma~\ref{lem:incomplemma}, either
   $X_1-x$ or $X_2-x$, say $X_1 - x$, would be a circuit-hyperplane of
   $M/x$. It follows that $X_1$ would be a circuit-hyperplane of $M$.
   This contradiction to the hypotheses of (\ref{thm:notchchar}.5)
   proves statement (i).  It follows that for each $x\in X_1\cap X_2$,
   the deletion $M\del x$ is a connected notch matroid, so by
   Lemma~\ref{lem:incomplemma}, either $X_1-x$ or $X_2-x$, say
   $X_1-x$, is a circuit-hyperplane of $M\del x$.  Since the circuit
   $X_1-x$ of $M\del x$ cannot contain parallel elements, statement
   (ii) follows.  By (\ref{thm:notchchar}.3) the minor $M|X_2/y\del x$
   is connected, so by part (b) of Corollary~\ref{cor:spanningvar}
   there is a spanning circuit $X'_2$ of $M|X_2$ that contains $y$.
   Lemma~\ref{lem:incomplemma} and the minimality of the excluded
   minor $M$ imply that $X_2$ is $X'_2\cup x$, so statement (iii)
   holds.  For statement (iv), note that if $|X_1|>|X_2|$ and $z\in
   X_1 - X_2$, then $M/z$, $X_1 - z$, and $\cl_{M/z}(X_2)$ contradict
   Lemma~\ref{lem:incomplemma}.  Statement (v) follows from part (iv)
   of Corollary~\ref{cor:finterval} since each of the notch matroids
   $M\del x$ and $M\del y$ has two circuit-hyperplanes.
   
   Now assume any two incomparable nontrivial connected flats are
   disjoint. We showed that the union of any two such flats is $E(M)$.
   Let $X_1,X_2$ be such flats.  It follows that all nonspanning
   circuits of $M$ span either $M|X_1$ or $M|X_2$, so $M$ is
   $T_n(M|X_1\oplus M|X_2)$; also, $M|X_1$ and $M|X_2$ are uniform
   matroids.  If $X_1$ is not a circuit and $x$ is in $X_1$, then
   $M\del x$ is a connected notch matroid in which $X_2$ is not a
   circuit-hyperplane, so $X_1 - x$ is a circuit-hyperplane of $M\del
   x$; it follows that $M|X_1$ is $U_{n-1,n+1}$. Assume that $X_1$ is
   a circuit, and so not a hyperplane of $M$; let $x$ be in $X_2$.
   Note that $X_1$ and $X_2-x$ are incomparable connected flats of the
   notch matroid $M/x$, which has no isthmuses. Since $X_2$ is not a
   circuit-hyperplane of $M$, it follows that $X_2-x$ cannot be a
   circuit-hyperplane of $M/x$.  Therefore by
   Lemma~\ref{lem:incomplemma}, $X_1$ is a circuit-hyperplane of
   $M/x$.  Thus, $M|X_1$ is $U_{n-2,n-1}$. In this manner, we see that
   there are, up to switching $X_1$ and $X_2$, three possibilities:
   $M|X_1$ and $M|X_2$ are both $U_{n-2,n-1}$; $M|X_1$ is
   $U_{n-2,n-1}$ and $M|X_2$ is $U_{n-1,n+1}$; both $M|X_1$ and
   $M|X_2$ are $U_{n-1,n+1}$. These possibilities give, respectively,
   $F_n$, $H_n$, and $G_n$.
\end{proof}

\begin{quote}
  {\bf (\ref{thm:notchchar}.6) } If relaxing some circuit-hyperplane
  $C$ of $M$ gives a generalized Catalan matroid $M'$, then $M$ is
  $A_n$.
\end{quote}

\begin{proof}[Proof of (\ref{thm:notchchar}.6).]
  We show that $M$ is $A_n$ by proving the following statements.
\begin{itemize}
  \item[(i)] There is a nonspanning circuit $C'\ne C$ of $M$ with $C\cap
    C'\ne\emptyset$.
  \end{itemize}
  Fix such a circuit $C'$ of least cardinality.
  \begin{itemize}
   \item[(ii)] There is at least one element
     $y$ in $E(M) - \bigl(C\cup\cl(C')\bigr)$.
   \item[(iii)] The ground set of $M$ is $C\cup C'\cup y$; also
     $|C\cap C'| = 1$.
  \item[(iv)] The circuit $C'$ is a hyperplane of $M$.
  \item[(v)] The only nonspanning circuits of $M$ are $C$ and
    $C'$.
  \end{itemize}
    
  Let the chain of proper nontrivial connected flats of $M'$ be
  $X_1\subset\cdots\subset X_k$. If $C\cap X_k$ were empty, then, by
  Corollary~\ref{cor:lpmauto}, there would be an automorphism of $M'$
  that maps $C$ to a final segment; by Lemma~\ref{lem:relax} we would
  get the contradiction that $M$ is a notch matroid.  Thus, $C\cap
  X_k$ is not empty, which gives statement (i).  Among all circuits
  that intersect $C$, choose $C'$ with smallest cardinality. The
  closure $\cl(C')$ is one of the connected flats $X_j$, and by the
  choice of $C'$, the basis $C$ of $M'$ is disjoint from $X_i$ for
  $i<j$. To prove statement (ii) we must show that $C$ does not
  contain the complement of $X_j$; if this were false, then by
  Corollary~\ref{cor:lpmauto} and Lemma~\ref{lem:relax} we would get,
  as before, that $M$ is a notch matroid.
  
  By Theorem~\ref{thm:notlpm}, $M|(C\cup C'\cup y)$ is not a lattice
  path matroid.  This observation and the minimality of $M$ prove the
  first part of statement (iii). The second part holds since if
  $|C\cap C'|\geq 2$ and $x\in C\cap C'$, then, by
  Theorem~\ref{thm:notlpm}, $M/x$ would not be a lattice path matroid.
  Let $C\cap C'$ be $x$.
  
  To prove statement (iv), first note that $M|\cl_M(C')$ is a uniform
  matroid since, by the choice of $C'$, any nonspanning circuit $Z$ of
  $M|\cl_M(C')$ would be disjoint from $C$, which gives the
  contradiction that the circuit $C'$ properly contains the circuit
  $Z$.  Since $M|\cl_M(C')$ is a uniform matroid that consists of $C'$
  and a subset of $C$, and since, by statement (iii), any circuit
  $C''\ne C$ with $|C''|=|C'|$ that intersects $C$ contains just one
  element of $C$, it follows that $C\cap \cl_M(C')$ is $x$, so $C'$ is
  closed.  If $C'$ is not a hyperplane of $M$, then there is an
  element $z$ in $C-\cl_M(C'\cup y)$, so $y$ is not in $\cl_M(C'\cup
  z)$.  However, for such a $z$, Theorem~\ref{thm:notlpm} applied to
  $M/z$, $\cl_{M/z}(C')$, $C-z$, and $y$ shows that $M/z$ is not in
  $\mathcal{L}$, contrary to $M$ being an excluded minor for
  $\mathcal{N}$.

  Since $C'$ is a circuit-hyperplane of $M$ and of the generalized
  Catalan matroid $M'$, it follows that $C'$ is the only nonspanning
  circuit of $M'$, so $C$ and $C'$ are the only nonspanning circuits
  of $M$, as needed to complete the proof.
\end{proof}
\end{proof}

\begin{figure}
\begin{center}
  \epsfxsize 2.6truein \epsffile{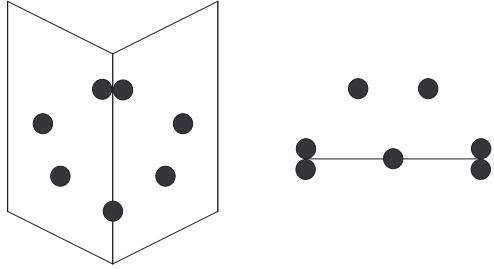}
\end{center}
\caption{Two more excluded minors for the class of lattice path
  matroids.}\label{otherex}
\end{figure}

Figure~\ref{otherex} shows two excluded minors for $\mathcal{L}$ that
are not among those given in Theorem~\ref{thm:exlexn}.  Presently we
do not know whether these two matroids complete the list of excluded
minors for the class of lattice path matroids.  

We close by noting that a lattice path matroid is graphic if and only
if it is the cycle matroid of an outerplanar graph in which each inner
face shares edges with at most two other inner faces.  One implication
follows since $\mathcal{W}_3$ and $C_{4,2}$ (i.e., the cycle matroids
of the two excluded minors, $K_4$ and $K_{2,3}$, for outerplanar
graphs) are excluded minors for lattice path matroids, as is
$B_{2,2}$, which is the cycle matroid of the graph formed by adding an
edge parallel to each edge of $K_3$.  The other implication follows
since by adding edges any graph of the stated type can be extended 
to a graph of this
type in which each face is bounded by at most three edges, and the
cycle matroids of such graphs, which are certain parallel connections
of $3$-point lines, are easily seen to be lattice path matroids.

\ 

\begin{center}
\textsc{Acknowledgements}
\end{center}

The authors thank Omer Gim\'enez for some useful observations related
to several parts of this paper.

\end{document}